\documentclass[11pt]{amsart}

\usepackage{amsmath}

\usepackage{amssymb}

\newcommand{\R}{{\mathbb R}}       
\newcommand{\DD}{{\mathcal D}}

\newcommand{\D}{{\Delta}}

\newcommand{\CC}{{\mathcal C}}

\newcommand{\dist}{{\rm dist}}
\newcommand{\ds}{\displaystyle }

\newcommand{\rf}[1]{{(\ref{#1})}}

\newcommand{\supp}{{\rm supp}}

\newcommand{\vphi}{{\varphi}}
\newcommand{\ve}{{\varepsilon}}
\newcommand{\vv}{{\vspace{2mm}}}

\newcommand{\wt}[1]{{\widetilde{#1}}}
\newcommand{\wh}[1]{{\widehat{#1}}}

\newcommand{\bmo}{{B\!M\!O}}
\newcommand{\rbmo}{{R\!B\!M\!O}}

\newcommand{\hbm}{{H^{1,\infty}_{atb}(\mu)}}

\newcommand{\QH}{{\wh{\wh{Q}^3 \hspace{-1mm}} }}

\newtheorem{theorem}{Theorem}[section]
\newtheorem{lemma}[theorem]{Lemma}

\newtheorem{propo}[theorem]{Proposition}
\newtheorem{claim}{Claim}

\theoremstyle{definition}
\newtheorem{definition}[theorem]{Definition}
\newtheorem{example}[theorem]{Example}

\theoremstyle{remark}
\newtheorem{rem}[theorem]{Remark}

\numberwithin{equation}{section}

\newcommand{\brem}{\begin{rem}}
\newcommand{\erem}{\end{rem}}

\newcommand{\bexam}{\begin{example}}
\newcommand{\eexam}{\end{example}}

\begin{document}

\title[Non doubling $H^1$ in terms of a maximal operator]
{Characterization of the atomic space $H^1$ for non doubling measures
in terms of a grand maximal operator}

\author[XAVIER TOLSA]{Xavier Tolsa}

\address{Department of Mathematics, Chalmers, 412 96 G\"oteborg, Sweden}

\email{xavier@math.chalmers.se}

\thanks{Supported by a postdoctoral grant from the European Comission for the
TMR Network ``Harmonic Analysis''. Also
partially supported by grants DGICYT PB96-1183 and CIRIT
1998-SGR00052 (Spain.}

\subjclass{Primary 42B20; Secondary 42B30}


\break

\begin{abstract}
Let $\mu$ be a Radon measure on $\R^d$, which may be non doubling.  The only
condition that $\mu$ must satisfy is the size condition $\mu(B(x,r))\leq C\,r^n$,
for some fixed $0<n\leq d$.
Recently, the author introduced spaces of type $\bmo(\mu)$ and $H^1(\mu)$ with
properties similar to ones of the classical spaces $\bmo$
and $H^1$ defined for doubling measures. These new spaces proved to be useful
to study the $L^p(\mu)$ boundedness of Calder\'on-Zygmund operators without
assuming doubling conditions.
In this paper a characterization of this new
atomic Hardy space $H^1(\mu)$ in terms of a maximal operator $M_\Phi$ is given.
It is shown that $f$ belongs to $H^1(\mu)$ if and only if $f\in
L^1(\mu)$, $\int f\, d\mu=0$ and $M_\Phi f \in L^1(\mu)$, as in the usual
doubling situation.
\end{abstract}

\keywords{BMO, atomic spaces, Hardy spaces, Calder\'on-Zygmund operators, non
doubling measures, maximal functions, grand maximal operator}

\maketitle


\section{Introduction}

The aim of this paper is to characterize the atomic Hardy space $\hbm$
introduced in \cite{Tolsa3} in terms of a grand maximal operator. Throughout all the
paper $\mu$ will be a (positive) Radon measure on $\R^d$ satisfying the growth
condition
\begin{equation}  \label{creix}
\mu(B(x,r))\leq C_0\,r^n  \qquad \mbox{for all $x\in \supp(\mu)$, $r>0$,}
\end{equation}
where $n$ is some fixed number with
$0<n\leq d$. We do {\em not} assume that $\mu$ is doubling
($\mu$ is said to be doubling if there exists some constant $C$
such that $\mu(B(x,2r)) \leq C\,\mu(B(x,r))$ for all $x\in \supp(\mu)$, $r>0$).

The doubling condition on $\mu$ is an essential assumption in most
results of classical Calder\'on-Zygmund theory. Nevertheless, recently it has
been shown that many results in this theory also hold without the doubling
assumption. For example, in \cite{Tolsa1} a $T(1)$ theorem and weak $(1,1)$
estimates for the Cauchy tranforms are obtained. For general Calder\'on-Zygmund
operators (CZO's) a $T(1)$ theorem in \cite{NTV1}, and weak
$(1,1)$ estimates and Cotlar's inequality in \cite{NTV2}
are proved. A $T(b)$ is also given in \cite{NTV*}.
For more results, see \cite{MMNO}, \cite{NTV4}, \cite{OP},
\cite{Tolsa2.5}, \cite{Tolsa3}, \cite{Tolsa5} and \cite{Verdera}, for example.

In \cite{Tolsa3} some variants of the classical spaces $\bmo(\mu)$ and
$H^1(\mu)$ are
introduced. These variants are denoted by $\rbmo(\mu)$ and $\hbm$ respectively.
There, it is shown that many of the properties fulfiled by
$\bmo(\mu)$ and $H^1(\mu)$ when $\mu$ is doubling are also satisfied by
$\rbmo(\mu)$
and $\hbm$ without assuming $\mu$ doubling. For example, the functions
from $\rbmo(\mu)$ fulfil a John-Nirenberg type inequality (see Section \ref{JJNN}
for the precise statement of this inequality), $\rbmo(\mu)$ is the dual of
$\hbm$, CZO's which are bounded in $L^2(\mu)$ are also bounded
from $\hbm$ into $L^1(\mu)$ and from $L^\infty(\mu)$ into $\rbmo(\mu)$ and, on
the other hand, any operator which is bounded from $\hbm$ into $L^1(\mu)$ and
from $L^\infty(\mu)$ into $\rbmo(\mu)$ is bounded in $L^p(\mu)$, $1<p<\infty$.

Let us remark that if $\mu$ is non doubling and one defines $\bmo(\mu)$ and
the atomic space $H^{1,\infty}_{at}(\mu)\equiv H^1(\mu)$
exactly as in the classical doubling situation (see
\cite{GR}, \cite{Journe} or \cite{Stein}, for instance), then these spaces
still fulfil some of the properties stated above \cite{MMNO}. However a
basic one fails: CZO's may be bounded in $L^2(\mu)$ but
not from $H^{1,\infty}_{at}(\mu)$ into $L^1(\mu)$ or
from $L^\infty(\mu)$ into $\bmo(\mu)$ (see \cite{Verdera} and \cite{MMNO}).
For this reason, if one wants to study the $L^p$-boundedness of CZO's,
the spaces $\bmo(\mu)$ and $H^{1,\infty}_{at}(\mu)$ are not appropiate.
This is the main reason for the introduction of $\rbmo(\mu)$ and $\hbm$ in
\cite{Tolsa3}.

\enlargethispage{10mm}
Before stating our main result, we need some notation and terminology.
By a cube $Q\subset \R^d$ we mean a closed cube centered at some point in $\supp(\mu)$ with
sides parallel to the
axes. Its side length is denoted by $\ell(Q)$ and its center by $z_Q$.
Given $\rho>0$, we denote by $\rho Q$ the cube concentric with $Q$ with
side length $\rho\,\ell(Q)$.
Recall that a function $f\in L^1_{loc}(\mu)$ belongs to the classical space
$H^{1,\infty}_{at}(\mu)$ if it can be written as $f=\sum_i \lambda_i\,a_i$,
where $\lambda_i\in\R$ are numbers such that $\sum_i |\lambda_i|<\infty$ and $a_i$
are functions called {\em atoms} such that
\begin{enumerate}
\item[1.] there exists some cube $Q_i$ such that $\supp(a_i)\subset Q_i$,
\item[2.] $\ds \int a_i\, d\mu = 0$,
\item[3.] $\|a_i\|_{L^\infty(\mu)} \leq \mu(Q_i)^{-1}.$
\end{enumerate}

In order to recall the precise definition of $\hbm$ we have to introduce
the coefficients $K_{Q,R}$.
Given two cubes $Q\subset R$, we set
$$K_{Q,R} = 1 + \int_{Q_R\setminus Q} \frac{1}{|x-z_Q|^n}\, d\mu(x),$$
where $Q_R$ is the smallest cube concentric with $Q$ containing $R$.

For a fixed $\rho>1$, a function $b\in L^1_{loc}(\mu)$ is called an {\em
atomic block} if
\begin{enumerate}
\item[1.] there exists some cube $R$ such that $\supp(b)\subset R$,
\item[2.] $\ds \int b\, d\mu = 0$,
\item[3.] there are functions $a_j$ supported on cubes $Q_j\subset R$ and
numbers $\lambda_j\in\R$ such that
$b=\sum_{j=1}^\infty \lambda_j a_j,$
and
$$\|a_j\|_{L^\infty(\mu)} \leq \left(\mu(\rho Q_j)\,K_{Q_j,R}\right)^{-1}.$$
\end{enumerate}
We denote
$$|b|_\hbm = \sum_j |\lambda_j|$$
(to be rigorous, we should think that $b$ is not only a function, but a
`structure'
formed by the function $b$, the cubes $R$ and $Q_j$, the functions $a_j$,
etc.). Then, we say that $f\in \hbm$ if there are atomic blocks $b_i$
such that
\begin{equation}  \label{sumatb}
f= \sum_{i=1}^\infty b_i,
\end{equation}
with $\sum_i |b_i|_\hbm<\infty$ (notice that this implies that the sum in
\rf{sumatb} converges in $L^1(\mu)$). The $\hbm$ norm of $f$ is
$$\|f\|_\hbm = \inf \sum_i |b_i|_\hbm,$$
where the infimum is taken over all the possible decompositions of $f$
in atomic blocks.

The definition of $\hbm$ does not depend on the constant $\rho>1$. The $\hbm$
norms for different choices of $\rho>1$ are equivalent. Nevertheless, for
definiteness, we will assume $\rho=2$ in the definition.

Compare the definitions of the spaces $H^{1,\infty}_{at}(\mu)$ and $\hbm$:
In $H^{1,\infty}_{at}(\mu)$ the cancellation condition 2 and the size condition
3 are imposed over the atoms $a_j$. On the other hand, in $\hbm$ the
cancellation condition 2 is imposed over the atomic blocks $b_i$, and the size
condition 3 is satisfied by the ``components'' $a_{i,j}$ of $b_i$ separately for
each $j$.
It is not difficult to check that $H^{1,\infty}_{at}(\mu) \equiv \hbm$ if
$\mu(B(x,r)) \approx r$ for all $x\in\supp(\mu),\,r>0$ (the notation
$A\approx B$
means that there exists some constant $C>0$ such that $C^{-1}\,A\leq B\leq
C\,A$, that is $A\lesssim B\lesssim A$). If the latter condition does not hold,
then $H^{1,\infty}_{at}(\mu)$ may be different from $\hbm$, even
when $\mu$ is doubling (see \cite{Tolsa3}).

Now we are going to introduce the ``grand'' maximal operator
$M_\Phi$, which is the main tool in our characterization  of $\hbm$.

\begin{definition} Given $f\in L^1_{loc}(\mu)$, we set
$$M_\Phi f(x) = \sup_{\vphi \sim x} \left| \int f\,\vphi\,d\mu\right|,$$
where the notation $\vphi \sim x$ means that $\vphi\in L^1(\mu)\cap
\CC^1(\R^d)$ and satisfies
\begin{enumerate}
\item $\|\vphi\|_{L^1(\mu)}\leq 1$,
\item $0\leq\vphi(y)\leq \dfrac{1}{|y-x|^n}$ for all $y\in \R^d$, and
\item $|\vphi'(y)|\leq \dfrac{1}{|y-x|^{n+1}}$ for all $y\in \R^d$.
\end{enumerate}
\end{definition}

In this paper we will prove the following result.

\begin{theorem}  \label{maxi}
A function $f$ belongs to $\hbm$ if and only if $f\in L^1(\mu)$,
$\int f\,d\mu=0$ and $M_\Phi f\in L^1(\mu)$. Moreover, in this case
$$\|f\|_{\hbm}\approx \|f\|_{L^1(\mu)} +  \|M_\Phi f\|_{L^1(\mu)}.$$
\end{theorem}


Theorem \ref{maxi} can be considered as a version for non doubling measures of
some results that are already known in more classical situations.
When $\mu$ is the Lebesgue measure on the real line, a characterization of
$H^{1,\infty}_{at}(\mu)$ such as the one
of Theorem \ref{maxi} was proved by Coifman \cite{Coifman}.
This result was extended to the Lebesgue measure on $\R^d$ by Latter
\cite{Latter}. Let us remark that in these cases, in the definition of $M_\Phi$,
for each $x$ it is enough to take the supremum over functions $\vphi_{x,r},
\,r>0$, of the form
$$\vphi_{x,r}(y) = \frac{1}{r^d}\,\psi\left(\frac{y-x}{r} \right),$$
where $0\not\equiv \psi\in {\mathcal S}$ is some fixed function.

If
\begin{equation}  \label{homnom}
\mu(B(x,r))\approx r^n \qquad \mbox{for all $x\in\supp(\mu),\,r>0$,}
\end{equation}
then $\supp(\mu)$ is a homogeneous space in the sense of \cite{CW}.
For general homogeneous spaces satisfying \rf{homnom}, Coifman,
Meyer and Weiss showed that there exists a description of
$H^{1,\infty}_{at}(\mu)$ in terms of a grand maximal operator 
(see \cite{CW} for this result and for the
detailed definition of homogeneous spaces). They observed
that a proof of this description by Carleson \cite{Carleson} using the duality
$H^{1,\infty}(\mu)$--$\bmo(\mu)$ in the case where $\mu$ is the Lebesgue
measure on $\R^n$ can be easily extended to the more general situation of
homogeneous spaces.

For a measure $\mu$ on $\R^d$ which is doubling but which may not satisfy
\rf{homnom}, Mac\'{\i}as and Segovia (\cite{MS1}, \cite{MS2}) obtained a
characterization of $H^{1,\infty}_{at}(\mu)$ by means of a grand maximal
operator too (see also \cite{Uchiyama}). They showed that if
$\mu$ is doubling, then taking a suitable quasimetric one can assume that
\rf{homnom} holds.
Their result applies not only to doubling measures on $\R^d$, but to more
general homogeneous spaces. On the other hand, since
$H^{1,\infty}_{at}(\mu)$ may be different from $\hbm$
if $\mu$ is a doubling measure on $\R^d$ which does not satisfy \rf{homnom},
the result of Mac\'{\i}as and Segovia (in the precise case that we are 
considering) cannot be derived as a particular instance of Theorem \ref{maxi}.

The absence of any regularity condition on $\mu$, apart from the size condition 
\rf{creix}, makes impossible to extend
the classical arguments to the present situation without major changes. We
will not consider any quasimetric on $\R^d$ different from the Euclidean
distance and we are not able to reduce our case to a situation
where \rf{homnom} holds.

Let us remark that the results of \cite{Coifman}, \cite{Latter}, \cite{MS1}
and \cite{MS2} concern not only the Hardy space $H^1$ but also
the Hardy spaces $H^p$, with $0<p<1$.
However, it is not possible to extend our proof of Theorem
\ref{maxi} to $0<p<1$
because we have obtained it by duality (following the same approach as Carleson
\cite{Carleson}).

\enlargethispage{10mm}

The paper is organized as follows. In Section \ref{secprelim} we deal with some
preliminary questions. In Section \ref{seceasy}
we show that the grand maximal operator $M_\Phi$ is bounded from $\hbm$ into
$L^1(\mu)$, which proves the ``only if'' part of Theorem \ref{maxi} (the
easy implication). 
In the remaining sections of the paper we prove the other
implication. In Section \ref{secdual} we explain how this
can be proved by duality. A suitable version for our purposes of 
John-Nirenberg inequality if obtained in Section \ref{JJNN}. 
In Section \ref{secdyad} some kind of dyadic cubes are constructed, 
and in the following section a suitable approximation 
of the identity adapted to the measure $\mu$ is obtained. 
Section \ref{seccore} contains a construction which is the core of the proof of the ``if'' part of
Theorem \ref{maxi}. 
Finally, Section \ref{secappen} is an Appendix where we prove a density result
which is necessary in the proof by duality of the ``if'' part of Theorem
\ref{maxi}.


\section{Preliminaries}   \label{secprelim}

The letter $C$ will be used for constants that may change from one occurrence to
another. Constants with subscripts, such as $C_1$, do not change in different
occurrences.

We will assume that the constant $C_0$ in \rf{creix} has been chosen big enough
so that for all the cubes $Q\subset \R^d$ we have
\begin{equation}\label{creix2}
\mu(Q) \leq C_0\, \ell(Q)^n.
\end{equation}
Given a function $f\in L^1_{loc}(\mu)$, we denote by $m_Qf$ the mean of $f$ over
$Q$ with respect to $\mu$, i.e. $m_Qf =
\frac{1}{\mu(Q)}\,\int_Q f\,d\mu$.

\begin{definition}
Given $\alpha>1$ and $\beta>\alpha^n$, we say that the cube $Q\subset \R^d$
is $(\alpha,\beta)$-doubling if $\mu(\alpha Q) \leq \beta\,\mu(Q)$.
\end{definition}

\begin{rem} \label{moltdob}
As shown in \cite{Tolsa3},
due to the fact that $\mu$ satisfies the growth condition \rf{creix}, there
are a lot ``big'' doubling cubes. To be precise, given any point
$x\in\supp(\mu)$ and $c>0$, there exists some
$(\alpha,\beta)$-doubling cube $Q$ centered at $x$ with $l(Q)\geq c$. This
follows easily from \rf{creix} and the fact that $\beta>\alpha^n$.

On the other hand, if $\beta>\alpha^d$, then for $\mu$-a.e. $x\in\R^d$ there
exists a sequence of $(\alpha,\beta)$-doubling cubes $\{Q_k\}_k$ centered
at $x$ with $\ell(Q_k)\to0$ as $k\to \infty$.
So there are a lot of ``small'' doubling cubes too.

For definiteness, if $\alpha$ and $\beta$ are not specified, by a doubling
cube we mean a $(2,2^{d+1})$-doubling cube.
\end{rem}

Now we are going to recall the definition of $\rbmo(\mu)$. In fact, in
Section 2 of \cite{Tolsa3} several equivalent definitions are given. Maybe  the
easiest one is the following. Let $f\in L^1_{loc}(\mu)$.
We say that $f\in\rbmo(\mu)$ if there exists some constant $C_1$ such that
for any {\em doubling} cube $Q$
\begin{equation} \label{def1}
\int_Q |f-m_Q f|\,d\mu \leq  C_1\, \mu(Q)
\end{equation}
and
\begin{equation} \label{def2}
|m_Q f-m_R f|\leq C_1 \,K_{Q,R}
\quad \mbox{for any two {\em doubling} cubes $Q\subset R$}.
\end{equation}
The best constant $C_1$ is the $\rbmo(\mu)$ norm of $f$, that we denote as
$\|f\|_*$.

Given any pair of constants $0<\alpha,\beta$, with $\beta>\alpha^n$,
if in the definition of $\rbmo(\mu)$ we ask \rf{def1} and \rf{def2} to hold for
$(\alpha,\beta)$-doubling cubes (instead of doubling cubes),
we will get the same space $\rbmo(\mu)$, with an equivalent norm \cite{Tolsa3}.
In fact, $\rbmo(\mu)$ can be defined also without talking about
doubling cubes: Given some fixed constant $\rho>1$, $f\in \rbmo(\mu)$ if and
only if there exists a collection of
numbers $\{f_Q\}_Q$ (i.e. for each cube $Q$ some number $f_Q$) and some constant
$C_2$ such that
$$\int_Q |f(x)-f_Q|\,d\mu(x) \leq C_2\,\mu(\rho Q)\quad \mbox{for any cube
$Q\subset \R^d$}$$
and,
$$|f_Q-f_R|\leq C_2\,K_{Q,R} \quad \mbox{for any two cubes $Q\subset R$}.$$
The best constant $C_2$ is comparable to the $\rbmo(\mu)$ norm of $f$ given by
\rf{def1} and \rf{def2}.

Recall that given two cubes $Q\subset R$, $Q_R$ stands for the smallest cube
concentric with $Q$ containing $R$. Without assuming $Q\subset R$,
we will denote by $Q_R$ the smallest cube concentric with $Q$ containing $Q$ and
$R$.

\begin{definition} \label{deltaqr}
Consider two cubes $Q,R\subset\R^d$ (we do {\em not} assume $Q\subset R$ or
$R\subset Q$). We denote
$$\delta(Q,R) = \max \left( \int_{Q_R\setminus Q} \frac{1}{|x-z_Q|^n}\, d\mu(x),
\,\,\int_{R_Q\setminus R} \frac{1}{|x-z_R|^n}\, d\mu(x) \right).$$
\end{definition}

Notice that $\ell(Q_R) \approx \ell(R_Q) \approx \ell(Q) + \ell(R) +
\dist(Q,R)$, and if $Q\subset R$, then $R_Q=R$ and $\ell(R)\leq \ell(Q_R) \leq
2\ell(R)$.

It is clear that if $Q\subset R$, then $K_{Q,R} = 1 + \delta(Q,R)$.
Quite often we will treat points $x\in\supp(\mu)$ as if they were cubes (with
$\ell(x)=0$). So for $x,y\in\supp(\mu)$ and some cube $Q$, the notations
$\delta(x,Q)$ and $\delta(x,y)$ make sense. In some way, they are particular
cases of Definition \ref{deltaqr}. Of course, it may happen
$\delta(x,Q) = \infty$ or $\delta(x,y)=\infty$.

In the following lemma we show that $\delta(\cdot,\cdot)$ satisfies some
very useful properties.

\begin{lemma}  \label{propdelta}
The following properties hold:
\begin{itemize}
\item[(a)] If $\ell(Q) \approx \ell(R)$ and $\dist(Q,R)\lesssim \ell(Q)$, then
$\delta(Q,R) \leq C$. In particular, $\delta(Q,\rho Q)\leq C_0\,2^n\,\rho^n$ for
$\rho>1$.

\item[(b)] Let $Q\subset R$ be concentric cubes such that there are no doubling
cubes of the form $2^k Q$, $k\geq0$, with $Q\subset 2^k Q\subset R$. Then,
$\delta(Q,R)\leq C_3$.

\item[(c)] If $Q\subset R$, then
$$\delta(Q,R) \leq C\,\left(1+\log\frac{\ell(R)}{\ell(Q)}\right).$$

\item[(d)] If $P\subset Q \subset R$, then
$$\bigl|\delta(P,R) - [ \delta(P,Q) + \delta(Q,R) ] \bigr| \leq \ve_0.$$
That is, with a different notation,
$\delta(P,R) = \delta(P,Q) + \delta(Q,R) \pm \ve_0$.
If $P$ and $Q$ are concentric, then
$\ve_0=0$: $\delta(P,R) = \delta(P,Q) + \delta(Q,R)$.

\item[(e)] For $P,Q,R\subset \R^d$,
$$\delta(P,R) \leq C_4\, + \delta(P,Q) + \delta(Q,R).$$
\end{itemize}
\end{lemma}

The constants that appear in (b), (c), (d) and (e) depend on $C_0,n,d$.
The constant $C$ in (a) depends, further, on the constants that are
implicit in the relations $\approx,\, \lesssim$.

Let us insist on the fact that a notation such as $a= b \pm \ve$ does not mean
any precise equality but the estimate $|a-b|\leq \ve$.

\begin{proof}
The estimates in (a) are immediate. The proof of (b) is also an easy estimate,
which can be found in \cite[Lemma 2.1]{Tolsa3}, for example.
The arguments for (c) are also quite standard. We leave the proof for
the reader.

Let us see that (d) holds.
If $P$ and $Q$ are concentric, the identity
$\delta(P,R) = \delta(P,Q) + \delta(Q,R)$ is a direct consequence of the
definition. In case $P$ and $Q$ are not concentric we have to make some
calculations:
\begin{eqnarray*}
\delta(P,R) & = & \delta(P,P_Q) + \int_{P_R\setminus P_Q}
\frac{1}{|y-z_P|^n}\,d\mu(y) \\
& = & \delta(P,Q) + \int_{P_R\setminus P_Q} \frac{1}{|y-z_P|^n}\,d\mu(y).
\end{eqnarray*}
So we must show that
$$S := \left|\int_{P_R\setminus P_Q} \frac{1}{|y-z_P|^n}\,d\mu(y) -
\delta(Q,R)\right| \leq C.$$
We set
\begin{eqnarray*}
S & \leq & \int_{P_Q\setminus Q} \frac{1}{|y-z_Q|^n}\,d\mu(y) +
\int_{P_R\D Q_R} \left( \frac{1}{|y-z_P|^n} + \frac{1}{|y-z_Q|^n} \right)\,
d\mu(y)\\
&&\mbox{} + \int_{\R^d\setminus P_Q} \left|\frac{1}{|y-z_P|^n} -
\frac{1}{|y-z_Q|^n} \right| \, d\mu(y) \\
& = & S_1 + S_2 + S_3.
\end{eqnarray*}
The integral $S_2$ is easily estimated above by some constant $C$, since
$|y-z_P|,\,|y-z_Q|\leq C\,\ell(R)$ for $y\in P_R\D Q_R$. An
analogous calculation yields $S_1\leq C$. For $S_3$ we have
$$S_3\leq C\,\int_{|y-z_Q|\geq \ell(Q)/2} \frac{|z_P-z_Q|}{|y-z_Q|^{n+1}}
\,d\mu(y)
\leq C\, \frac{|z_P-z_Q|}{\ell(Q)} \leq C,$$
and we are done with (d).

We leave the proof of (e) for the reader too.
\end{proof}

Notice that if we set $D(Q,R) = 1 + \delta(Q,R)$ for $Q\neq R$ and
$D(Q,Q) = 0$, then $D(\cdot,\cdot)$ is a quasidistance on the set of cubes,
by (e) in the preceding lemma.

From (a) and the fact that $Q_R$ and $R_Q$ have comparable sizes and
$Q_R\cap R_Q\neq\varnothing$, we get that $Q_R$ and $R_Q$ are close in the
quasimetric $D(\cdot,\cdot)$.
Also,  if we denote by $\wt{Q}$ the smallest doubling cube of the form $2^kQ$,
$k\geq0$, by (b) we know that $\wt{Q}$ is not far from $Q$ (using again the
quasidistance $D$). So $Q$ and $\wt{Q}$ may have very different sizes, but we
still have $D(Q,\wt{Q})\leq C$.

\vv
In Remark \ref{moltdob} we have explained that there a lot of big and small
doubling cubes. In the following lemma we state a more precise result about the
existence of small doubling cubes in terms of $\delta(\cdot,\cdot)$.

\begin{lemma}\label{precisio}
There exists some (big) constant $\eta>0$ depending only on $C_0$, $n$ and $d$
such that if $R_0$ is some cube centered at some point of $\supp(\mu)$ and
$\alpha>\eta$, then
for each $x\in R_0\cap \supp(\mu)$ such that $\delta(x,2R_0) >\alpha$ there
exists some doubling cube $Q\subset 2R_0$ centered at $x$ satisfying
\begin{equation}  \label{err1}
|\delta(Q,2R_0) - \alpha|\leq \ve_1,
\end{equation}
where $\ve_1$ depends only on $C_0$, $n$ and $d$ (but not on $\alpha$).
\end{lemma}

\begin{proof}
Let $Q_1$ be the biggest cube centered at $x$ with side length
$2^{-k}\,\ell(R_0)$,
$k\geq1$, such that $\delta(Q_1,2R_0)\geq\alpha$. Then,
$\delta(2Q_1,2R_0)<\alpha$. Otherwise, $k=1$ and since $\ell(Q_1)= \ell(R_0)/2$
and $\ell(Q_{1,R_0})\leq 4\,\ell(R_0)$ we get
$$\delta(Q_1,2R_0) \leq \int_{\ell(Q_1)/2<|y-x|,\, y\in Q_{1,R_0}}
\frac{1}{|y-x|^n}\,d\mu(y)
\leq \frac{C_0\,8^n \,\ell(R_0)^n}{\ell(Q_1)^n} = C_0\,16^n,$$
which contradicts the choice of $Q_1$, assuming $\eta>C_0\,16^n$.

Now we have
$\delta(Q_1,2R_0) \leq \alpha + \delta(Q_1,2Q_1) \leq  \alpha + C_0\,16^n.$
Thus
$$|\delta(Q_1,2R_0)- \alpha|\leq C_0\,16^n.$$

Let $Q$ be the smaller doubling cube of the form $2^k\,Q_1$, $k\geq 0$.
Then $\delta(Q_1,Q)\leq C_3$.
Also, $\ell(Q)\leq \ell(R_0)$. Otherwise, $R_0\subset 3Q$ and
$$\delta(Q_1,2R_0)\leq \delta(Q_1,3Q) = \delta(Q_1,Q) + \delta(Q,3Q) \leq C_3
+ 6^n\,C_0.$$
This is not possible if we assume $\eta>C_3 + 6^n\,C_0$.

Now $Q$ satisfies the required properties, since it is doubling, it is contained
in $2R_0$, and
\begin{eqnarray*}
|\delta(Q,2R_0) - \alpha| & \leq &
|\delta(Q,2R_0) - \delta(Q_1,2R_0)| + |\delta(Q_1,2R_0) - \alpha|\\
& \leq & \delta(Q,Q_1)+  C_0\,16^n \leq C_3 + C_0\,16^n =:\ve_1.
\end{eqnarray*}
\end{proof}

As in (d) of Lemma \ref{propdelta},
instead of \rf{err1}, often we will write $\delta(Q,2R_0) = \alpha
\pm\ve_1$.

Notice that by (e) and (a) of Lemma \ref{propdelta}, we get
\begin{eqnarray*}
|\delta(Q,R_0) - \alpha| & \leq & |\delta(Q,2R_0) - \alpha| +
|\delta(Q,2R_0) - \delta(Q,R_0)| \\
& \leq & \ve_1 + \delta(R_0,2R_0) +C_4 \\
& \leq & \ve_1 + C +C_4 :=\ve_1'.
\end{eqnarray*}
However we prefer the estimate \rf{err1},
because we have $Q\subset 2R_0$ but $Q\not \subset R_0$, in general. So the
cube $2R_0$, in some sense, is a more appropriate reference.

Results analogous to the ones in Lemma \ref{precisio} can be stated about the
existence of cubes $Q$ centered at some point $x\in R_0$ with $Q\supset R_0$,
but since we will not need this fact below, we will not show any precise result
of this kind.

\vv
If $Q\subset R$ are doubling cubes and $f\in\rbmo(\mu)$, then
$|m_Qf - m_Rf|\leq (1+\delta(Q,R))\,\|f\|_*.$ Without assuming $Q\subset R$, we
have a similar result:

\begin{propo}
Let $Q,R\subset \R^d$ be doubling cubes. If $f\in \rbmo(\mu)$, then
$$|m_Qf - m_Rf|\leq (C+2\,\delta(Q,R))\,\|f\|_*.$$
\end{propo}

\begin{proof}
Suppose, for example, $\ell(R_Q)\geq \ell(Q_R)$. Then, $Q_R\subset 3R_Q$.

Let $\wt{3R_Q}$ be the smallest doubling cube of the form
$2^k\,3R_Q$, $k\geq0$. We have
$$\delta(R,\wt{3R_Q}) = \delta(R,R_Q) +
\delta(R_Q,\wt{3R_Q})  \leq \delta(R,Q) + C.$$
Thus
\begin{equation}  \label{acott1}
|m_Rf - m_{\wt{3R_Q}}f| \leq (1 + C + \delta(R,Q))\, \|f\|_*.
\end{equation}
We also have
$$\delta(Q,\wt{3R_Q}) \leq  C + \delta(Q,3R_Q) + \delta(3R_Q,\wt{3R_Q})
\leq C + \delta(Q,Q_R) +  \delta(Q_R,3R_Q).$$
Since $Q_R$ and $R_Q$ have comparable sizes, $\delta(Q_R,3R_Q) \leq C$, and so
$$\delta(Q,\wt{3R_Q}) \leq C + \delta(Q,R).$$
Therefore,
\begin{equation}  \label{acott2}
|m_Qf - m_{\wt{3R_Q}}f| \leq (1 + C + \delta(Q,R))\, \|f\|_*.
\end{equation}
By \rf{acott1} and \rf{acott2}, the proposition follows.
\end{proof}


\section{The easy implication of Theorem \ref{maxi}}  \label{seceasy}

In this section we will prove the ``only if'' part of Theorem
\ref{maxi}.

\begin{lemma} \label{only}
The operator $M_\Phi$ is bounded from $\hbm$ into $L^1(\mu)$.
\end{lemma}

\begin{proof}
Let $b=\sum_i \lambda_i\, a_i$ be an atomic block supported on some cube $R$,
with $\lambda_i\in\R$, where $a_i$ are functions
supported on cubes $Q_i\subset R$ such that $\|a_i\|_{\infty} \leq
((1 + \delta(Q_i,R))\, \mu(2Q_i))^{-1}$.
We will show that $\|M_\Phi b\|_{L^1(\mu)} \leq C\,\sum_i|\lambda_i|$.

First we will estimate the integral
$\int_{\R^d\setminus 2R} M_\Phi b\,\, d\mu.$
For $x\in \R^d\setminus 2R$ and $\vphi\sim x$, since $\int b\,d\mu=0$, we
have
\begin{eqnarray}  \label{tyy1}
\left|\int b\,\vphi\,d\mu\right| & = & \left|\int
b(y)\,(\vphi(y)-\vphi(z_R))\,d\mu(y)\right| \nonumber \\
& \leq & C\,\int |b(y)|\,\frac{\ell(R)}{|x-z_R|^{n+1}}\,d\mu(y).
\end{eqnarray}
Thus
\begin{eqnarray} \label{tyy2}
\int_{\R^d\setminus 2R} M_\Phi b\,\, d\mu
& \leq & C\,\|b\|_{L^1(\mu)} \int_{\R^d\setminus 2R}
\frac{\ell(R)}{|x-z_R|^{n+1}}\,d\mu(x) \nonumber \\
& \leq & C\,\|b\|_{L^1(\mu)} \leq C\,
\sum_i|\lambda_i|.
\end{eqnarray}

Now we will show that
\begin{equation} \label{eq1}
\int_{2R} M_\Phi a_i\, d\mu \leq C,
\end{equation}
and we will be done.
If $x\in 2Q_i$ and $\vphi\sim x$, then
$$\left|\int a_i\,\vphi\,d\mu\right| \leq C\,\|a_i\|_{L^\infty(\mu)}\,
\|\vphi\|_{L^1(\mu)} \leq C\,\|a_i\|_{L^\infty(\mu)}.$$
So
$$\int_{2Q_i} M_\Phi a_i\, d\mu \leq C\,\|a_i\|_{L^\infty(\mu)}\,\mu(2Q_i)
\leq C.$$
For $x\in 2R\setminus 2Q_i$ and $\vphi\sim x$, we have
$$\left|\int a_i\,\vphi\,d\mu\right| \leq
C\,\|a_i\|_{L^1(\mu)}\,\frac{1}{|x-z_{Q_i}|^n}.$$
Therefore,
\begin{eqnarray}  \label{tyy3}
\int_{2R\setminus 2Q_i} M_\Phi a_i\, d\mu &\leq &
C\,\|a_i\|_{L^1(\mu)}\, \int_{2R\setminus 2Q_i}\frac{1}{|x-z_{Q_i}|^n}\,
d\mu(x) \nonumber \\
& \leq & C\,\|a_i\|_{L^1(\mu)}\, (1+ \delta(Q_i,R)) \leq C,
\end{eqnarray}
and \rf{eq1} follows.
\end{proof}


\section{An approach by duality for the other implication} \label{secdual}

We have to show that if $f\in L^1(\mu)$,
$\int f\,d\mu=0$ and $M_\Phi f\in L^1(\mu)$, then $f\in\hbm$.
We will obtain this result by duality, following the ideas of Carleson
\cite{Carleson}.
So we will prove

\begin{lemma}[\bf Main Lemma]
Let $f\in \rbmo(\mu)$ with compact support and
$\int f\,d\mu=0$. There exist functions
$h_m\in L^\infty(\mu)$, $m\geq 0$, such that
\begin{equation}  \label{main}
f(x) = h_0(x) + \sum_{m=1}^\infty \int \vphi_{y,m}(x)\, h_m(y)\,d\mu(y),
\end{equation}
with convergence in $L^1(\mu)$ where, for each $m\geq1$, 
$\vphi_{y,m}\sim y$, and
\begin{equation}  \label{main2}
\sum_{m=0}^\infty |h_m| \leq C\,\|f\|_*.
\end{equation}
\end{lemma}

Let us see that from this lemma the ``if'' part of Theorem
\rf{maxi} follows. Consider $f\in L^1(\mu)$ such that $\int f\,d\mu=0$ and
$M_\Phi f\in L^1(\mu)$. Assume first that $f\in L^\infty(\mu)$ and has compact
support. In this case, $f\in \hbm$ and so we only have to estimate the norm of
$f$.

Since $\rbmo(\mu)$ is the dual of $\hbm$ \cite{Tolsa3},
given $f\in L^1(\mu)$, by the Hahn-Banach theorem we have
$$\|f\|_\hbm = \sup_{\|g\|_*\leq1} |\langle f,\,g\rangle|.$$
Since $\int f\,d\mu=0$, we can assume that $g$ has compact support and
$\int g\,d\mu=0$. Then, applying the Main Lemma to $g$ we get
$$|\langle f,\,g\rangle| \leq \left|\int f\,h_0\,d\mu\right|
+ \left|\sum_{m=1}^\infty \iint \vphi_{y,m}(x)\,h_m(y)\,f(x)\,d\mu(x)\,d\mu(y)
\right|.$$
Since $\int |\vphi_{y,m}(x)\,f(x)|\,d\mu(x)\leq M_\Phi f(y)$, we have
\begin{eqnarray*}
|\langle f,\,g\rangle| & \leq & \|f\|_{L^1(\mu)}\, \|h_0\|_{L^\infty(\mu)}
+ \sum_{m=1}^\infty \int M_\Phi f(y)\,|h_m(y)|\,d\mu(y) \\
&\leq & \|f\|_{L^1(\mu)}\, \|h_0\|_{L^\infty(\mu)} + \|M_\Phi f\|_{L^1(\mu)}\,
\left\|\sum_{m=1}^\infty |h_m|\right\|_{L^\infty(\mu)} \\
& \leq & C\,\left(\|f\|_{L^1(\mu)}+ \|M_\Phi f\|_{L^1(\mu)}\right)\, \|g\|_*.
\end{eqnarray*}
That is, $\|f\|_\hbm \leq C\,\left(\|f\|_{L^1(\mu)}+
\|M_\Phi f\|_{L^1(\mu)}\right)$.

In the general case where we don't know a priori that $f\in\hbm$, we can
consider a sequence of functions $f_n$ bounded with compact support  such that
$\int f_n\,d\mu=0$, $f_n\to f$ in $L^1(\mu)$ and $\|M_\Phi
(f-f_n)\|_{L^1(\mu)}\to0$, and then we apply the usual arguments. The existence
of such a sequence is showed in Lemma \ref{final}, in the Appendix.

The rest of the paper, with the exception of the Appendix,
is devoted to the proof of the Main Lemma.


\section{The inequality of John-Nirenberg} \label{JJNN}

In \cite{Tolsa3} it is shown that the functions of the space $\rbmo(\mu)$
satisfy a John-Nirenberg type inequality. Let us state the precise result.

\begin{theorem} \label{JN1}
Let $Q\subset \R^d$ be a doubling cube. If $f\in\rbmo(\mu)$, then
$$\mu\{x\in Q:\, |f-m_Q f|>\lambda\} \leq C_5\,\mu(Q)\, \exp\left(
\frac{-C_6\,\lambda}{\|f\|_*}\right), \qquad \lambda>0,$$
where $C_5, C_6>0$ are constants that only depend on $C_0,\,n,\,d$.
\end{theorem}

In the proof of the Main Lemma we will need a version of the above inequality
which appears to be stronger (although it is equivalent). In this section we
will state and prove this new version of John-Nirenberg inequality.

\begin{definition}  \label{zq}
Given a doubling cube $Q$, we denote by $Z(Q,\lambda)$ the set of
points $x\in Q$
such that any doubling cube $P$ with $x\in P$ and $\ell(P)\leq \ell(Q)/4$ satisfies
$|m_P f - m_{Q}f|\leq \lambda.$
\end{definition}

In other other words, $Q\setminus Z(Q,\lambda)$ is the subset of $Q$
such that for some doubling cube $P$ with $x\in P$ and $\ell(P)\leq \ell(Q)/4$ we have
$$|m_P f - m_{Q}f|> \lambda.$$

\begin{propo}  \label{JN2}
Let $Q\subset \R^d$ be a doubling cube. If $f\in\rbmo(\mu)$, then
$$\mu(Q\setminus Z(Q,\lambda)) \leq C_5'\,\mu(Q)\, \exp\left(
\frac{-C_6'\,\lambda}{\|f\|_*}\right), \qquad \lambda>0.$$
where $C_5', C_6'>0$ are constants that  only depend on $C_0,\,n,\,d$.
\end{propo}

\begin{proof} The arguments are quite standard.
For any $x\in Q\setminus Z(Q,\lambda)$ there exists some cube $P_x$ which contains $x$,
with $\ell(P_x)\leq \ell(Q)/4$ and such that $|m_{P_x} f - m_{Q}f|> \lambda$.
Then by Besicovich's Covering Theorem, there are points
$x_i\in Q\setminus Z(Q,\lambda)$ such that
$$Q\setminus Z(Q,\lambda)\subset \bigcup_i 2P_i,$$
and so that the cubes $2P_i$, $i=1,2,\ldots$, form an almost disjoint family.
Observe that the Covering Theorem of Besicovich cannot be applied to the cubes
$P_x$ (they are non centered), however we have applied it to the cubes $2P_x$,
which are non centered too, but fulfil the condition
$$x\in \tfrac{1}{2} 2P_x.$$
That is, the point $x$ is ``far'' from the boundary of $2P_x$. Under this
condition, Besicovich's Covering Theorem also holds.

Since, for each $i$, $\ell(P_i)\leq \ell(Q)/4$ and $P_i\cap Q\neq\varnothing$, it is
easlily seen that $2P_i \subset \frac{7}{4}Q$. Then,
\begin{eqnarray*}
\mu(Q\setminus Z(Q,\lambda)) & \leq & \sum_i \mu(2P_i) \\
& \leq & \sum_i \int_{P_i} \exp\left(|f(x)-m_Q f|\,k\right)\,
\exp(-\lambda\,k)\, d\mu(x) \\
& \leq & C\,\int_{\frac{7}{4}Q} \exp\left(|f(x)-m_Q f|\,k\right)\,
 \exp(-\lambda\,k)\,d\mu(x),
\end{eqnarray*}
where $k$ is some constant that will be fixed below. Now, we have
\begin{eqnarray*}
\exp\left(|f(x)-m_Q f|\,k\right) & \leq &
\exp\left(|f(x)-m_{\frac{7}{4}Q} f|\,k\right)
\exp\left(|m_{\frac{7}{4}Q} f-m_Q f|\,k\right) \\
& \leq &
\exp\left(|f(x)-m_{\frac{7}{4}Q} f|\,k\right)
\exp\left(C\,\|f\|_*\,k\right).
\end{eqnarray*}
The last inequality follows from
$|m_{\frac{7}{4}Q} f-m_Q f| \leq C\,\|f\|_*$
(notice that the cube $\frac{7}{4}Q$ is $(\frac{8}{7},2^{d+1})$-doubling).

Therefore, by Theorem \ref{JN1} (which also holds for cubes that are
$(\frac{8}{7},2^{d+1})$-doubling instead of $(2,2^{d+1})$-doub\-ling, with
constants $\wt{C}_1$ and $\wt{C}_2$ instead of $C_1$ and $C_2$) we have
\begin{eqnarray*}
\lefteqn{
\mu(Q\setminus Z(Q,\lambda))} \\
 & \leq & C\,\exp(-\lambda\,k) \,
\exp\left(C\,\|f\|_*\,k\right)\,
\int_{\frac{7}{4}Q} \exp\left(|f(x)-m_{\tfrac{7}{4} Q} f|\,k\right)\,d\mu(x) \\
& = & C\,\exp(-\lambda\,k) \,
\exp\left(C\,\|f\|_*\,k\right)\\
&&\times \,\int_0^\infty \mu\left\{x\in \tfrac{7}{4}Q :\,
\exp\left(|f(x)-m_{\tfrac{7}{4} Q} f|\,k\right) > t\right\}\, dt \\
& \leq &
C\,\mu(\tfrac{7}{4}Q)\,
 \exp(-\lambda\,k) \, \exp\left(C\,\|f\|_*\,k\right)
\int_0^\infty \wt{C}_1\,\exp\left(\frac{-\wt{C}_2\,\log t}{k\,\|f\|_*}\right)
\, dt.
\end{eqnarray*}
So if we choose $k:=\wt{C}_2/2\|f\|_*$, we get
$$\mu(Q\setminus Z(Q,\lambda)) \leq C\,\mu(\tfrac{7}{4}Q) \,
\exp\left(\frac{-\wt{C}_2\,\lambda}{2\|f\|_*}\right)
\leq C\,\mu(Q)\,
\exp\left(\frac{-\wt{C}_2\,\lambda}{2\|f\|_*}\right).$$
\end{proof}


\section{The ``dyadic'' cubes}  \label{secdyad}

In \cite{Carleson}, Carleson proves a result analogous to the one stated in the
Main Lemma for $\mu$ being the
Lebesgue measure on $\R^d$. He uses dyadic cubes of
side length $2^{-mA}$, where $A$ is some big positive integer. In our proof,
we will also consider some cubes which will play the role of the dyadic cubes
with side length $2^{-mA}$ of Carleson. In this section we will introduce
these new ``dyadic'' cubes and we will show some of the properties that they
satisfy and that will be needed in the proof of the Main Lemma.

As in \cite{Carleson}, we will take some big positive integer $A$
whose precise value will be fixed after knowing or choosing several additional
constants. In particular, we assume that $A$ is much bigger than the constants
$\ve_0,\,\ve_1$ and $\eta$ of Section \ref{secprelim}.

\begin{definition}  \label{defcubs}
Suppose that the support of the function $f$ of the Main Lemma is contained in
a doubling cube $R_0$. Let $m\geq1$ be some fixed integer and
$x\in \supp(\mu)\cap R_0$. If $\delta(x,2R_0)> m\,A$, we denote by $Q_{x,m}$ a
doubling cube (with $Q_{x,m}>0$) such that
\begin{equation} \label{estimprec}
|\delta(Q_{x,m},2R_0) - m\,A|\leq \ve_1.
\end{equation}
Also, $\DD_m^\prime =\{Q_{i,m}\}_{i\in I^\prime_m}$, is a subfamily with
finite overlap of the
cubes $Q_{x,m}$, such that each cube $Q_{i,m}\equiv Q_{y_i,m}$ is centered at some point $y_i\in
\supp(\mu)\cap R_0$ with $\delta(y_i,2R_0)> m\,A,$ and
$$\{x\in \supp(\mu)\cap R_0:\,\delta(x,2R_0)> m\,A\} \subset \bigcup_{i\in
I^\prime_m} Q_{i,m}$$
(this family exists because of Besicovich's Covering Theorem).

If $\delta(x,2R_0)\leq m\,A,$ we set $Q_{x,m}=\{x\}$.
We denote by $\DD^{\prime\prime}_m$ the family of cubes $Q_{x,m}\equiv \{x\}$
such that $\delta(x,2R_0)\leq m\,A$ and $x\not \in \bigcup_{i\in
I^\prime_m} Q_{i,m}$. We set $\DD_m =\DD^\prime_m\cup\DD^{\prime\prime}_m.$

The cubes $Q_{x,m}$, $x\in\supp(\mu)\cap R_0$ (not necessarily from the family
$\DD_m$) are called cubes of the $m$-th generation.
\end{definition}

Obviously, the whole family of cubes in $\DD_m$ has also finite overlap.
Notice that if $x$ is a point in $\supp(\mu)$ such that
$\delta(x,2R_0) =\infty$, then $\ell(Q_{x,m})>0$ for all $m\geq1$.
Otherwise, there exists some $m_0$ such that $\ell(Q_{x,m})=0$ for all
$m\geq m_0$.

It is easily seen that if $A$ is big enough, then  $\ell(Q_{x,m+1})\leq
\ell(Q_{x,m})/10$
(a more precise version of this result will be proved in Lemma \ref{mides}
below). So $\ell(Q_{x,m}) \to 0$ as $m\to \infty$.

If $A$ is much bigger than $\ve_1$ and $Q_{x,m}\neq\{x\}$,
then $\delta(Q_{x,m},2R_0)\approx mA$.
However, the estimate \rf{estimprec} is much sharper.
This will very useful in our construction.

\begin{lemma} \label{idemgen}
Assume that $P$ and $Q$ are cubes contained in $2R_0$ whose centers are in
$R_0$. Let $S$ be a cube such that $P,Q\subset S\subset 2R_0$.
\begin{itemize}
\item[(a)] If
$|\delta(P,2R_0) - \delta(Q,2R_0)| \leq \beta,$
then
$$|\delta(P,S) - \delta(Q,S)| \leq \beta + 2\ve_0.$$

\item[(b)] If
$|\delta(P,S) - \delta(Q,S)|\leq \beta,$
then
$$|\delta(P,2R_0) - \delta(Q,2R_0)|\leq \beta + 2\ve_0.$$
\end{itemize}
\end{lemma}

In particular, this lemma can be applied to cubes $P$ and $Q$ belonging to
the same generation $m$, with $\beta=2\ve_1$ (assuming $\ell(P),\ell(Q)\neq0$).

\begin{proof}
Both statements are a straightforward consequence of (d) in Lemma \ref{propdelta}, since
$$\delta(P,2R_0) = \delta(P,S) + \delta(S,2R_0) \pm \ve_0$$
and
$$\delta(Q,2R_0) = \delta(Q,S) + \delta(S,2R_0) \pm \ve_0.$$
\end{proof}

The constants $\ve_0$ and $\ve_1$ should be understood as upper
bounds for some ``errors'' and deviations of our construction from the classical
dyadic lattice.

\vv
We will need the following result too.

\begin{lemma}  \label{mides}
Assume that $A$ is big enough. There exists some $\gamma>0$ such that
if $Q_{x,m} \cap Q_{y,m+1} \neq \varnothing$,
$x,y\in\supp(\mu)$, then $\ell(Q_{y,m+1}) \leq 2^{-\gamma\,A}\,\ell(Q_{x,m})$.
\end{lemma}

\begin{proof} We can assume $Q_{y,m+1}\neq \{y\}$. Let $B>1$ be some fixed
constant.
If $\ell(Q_{y,m+1}) > B^{-1}\,\ell(Q_{x,m})$, then $Q_{x,m}\subset 3B\,Q_{y,m+1}$.
So, if $R_x$ is a cube centered at $x$ with side length $6B\,\ell(Q_{y,m+1})$,
we have $Q_{x,m},Q_{y,m+1} \subset R_x$.

By (c) of Lemma \ref{propdelta} we get
$$\delta(Q_{y,m+1}, R_x) \leq C \,\left(1 +
\log\left(\frac{\ell(R_x)}{\ell(Q_{y,m+1})}\right) \right) \leq
C \,(1 + \log B) .$$
Since
$$\delta(Q_{y,m+1}, 2R_0) = \delta(Q_{y,m+1}, R_x) + \delta(R_x,2R_0)\pm\ve_0,$$
if we set $B= 2^{\gamma\,A}$, we obtain
$$\delta(R_x,2R_0) > (m+1)\,A - \ve_1 - \ve_0 - C (1+ \gamma\,A\, \log2).$$
Then for $\gamma$ small enough we have
$$\delta(R_x,2R_0) > (m+1)\,A - \ve_1 - \ve_0 - C - \frac{1}{2}\,A
> mA+\ve_1.$$
This implies $\delta(Q_{x,m},2R_0)> mA+\ve_1$, which is not possible.
\end{proof}

As a consequence, we obtain

\begin{lemma}  \label{mides2}
Assume that $A$ is big enough.
If $x,y\in\supp(\mu)$ are such that $Q_{x,m} \cap Q_{y,m+k} \neq \varnothing$
(with $k\geq1$),
then $\ell(Q_{y,m+k}) \leq 2^{-\gamma\,A\,k}\,\ell(Q_{x,m})$.
\end{lemma}

\begin{proof}
By the previous lemma, $\ell(Q_{y,j+1}) \leq 2^{-\gamma\,A}\,\ell(Q_{y,j})$ and
$\ell(Q_{y,m+1}) \leq 2^{-\gamma\,A}\,\ell(Q_{x,m})$. This gives
$\ell(Q_{y,m+1}) \leq 2^{-\gamma\,A\,k}\,\ell(Q_{x,m})$.
\end{proof}


\section{An approximation of the identity}  \label{secapide}

The proof of the Main Lemma will be constructive. At the level of cubes of
generation $m$ we will construct a function $h_m$ yielding the
``potential''
$$U_m(x) = \int \vphi_{y,m}(x)\,h_m(y)\,d\mu(y)$$
(to be precise, instead of one function $h_m$, for each $m$ we will have $N$
functions $h_m^1,\ldots, h_m^N$, but this is a rather technical detail that
we can skip now).
The potentials $U_m$ will compensate the large values of $f$ at the scale of
cubes of the generation $m$.
So the arguments will be similar to the ones of \cite{Carleson}.

However, in our situation several problems arise, in general, because of the
absence of any kind of regularity in the measure $\mu$ (except the
growth condition \rf{creix}). For example, in
\cite{Carleson} the potentials $U_m$ are convolutions with approximations of the
identity: $U_m = \vphi_m \ast h_m$.
Using the previous notation, we have
$$\vphi_{y,m}(x) = \vphi_m(y-x) = 2^{mAn}\,\vphi(2^{mA}(y-x)).$$
This is not our case. The measure $\mu$ is not invariant by translations
and we don't know how it behaves under dilations
(notice that if $\mu$ were doubling, we would have some information, at
least, about the behaviour under dilations).
We need to use
functions $\vphi_{y,m}$ such that $\|\vphi_{y,m}\|_{L^1(\mu)}=1$
(or at least equal to some value close to $1$).
So $\vphi_{y',m}$ cannot be obtained as a translation of
$\vphi_{y,m}$ for $y'\neq y$, neither as a dilation of $\vphi_{y',k}$, $k\neq m$.
In this section we will show how these problems can be overcome.

\vv
We denote
$$\sigma:= 10 \ve_0 + 10\ve_1 + 12^{n+1} C_0.$$

We introduce two new constants $\alpha_1, \alpha_2>0$ whose precise value will
be fixed below. For the moment, let us say that $\ve_0, \ve_1, C_0,
\sigma \ll \alpha_1\ll\alpha_2\ll A$.

\begin{definition} Let $y\in \supp(\mu)$.
We denote by $Q_{y,m}^1$, $\wh{Q}_{y,m}^1$, $Q_{y,m}^2$, $\wh{Q}_{y,m}^2$,
$Q_{y,m}^3$ some
doubling cubes (with positive side length) centered at $y$ such that
\begin{equation}  \label{igu}
\begin{split}
& \delta(Q_{y,m},2R_0) = m\,A \pm \ve_1,\\
& \delta(Q_{y,m}^1,2R_0) = m\,A - \alpha_1 \pm \ve_1,\\
& \delta(\wh{Q}_{y,m}^1,2R_0) = m\,A - \alpha_1 -\sigma \pm \ve_1,\\
& \delta(Q_{y,m}^2,2R_0) = m\,A - \alpha_1 - \alpha_2 \pm \ve_1,\\
& \delta(\wh{Q}_{y,m}^2,2R_0) = m\,A - \alpha_1 - \alpha_2 -\sigma \pm \ve_1 ,\\
& \delta(Q_{y,m}^3,2R_0) = m\,A - \alpha_1 - \alpha_2 -2\,\sigma \pm \ve_1
\end{split}
\end{equation}
By Lemma \ref{precisio} we know that if $\delta(y,2R_0)>m\,A$,
then all the cubes $Q_{y,m}^1$, $\wh{Q}_{y,m}^1$, $Q_{y,m}^2$,
$\wh{Q}_{y,m}^2$, $Q_{y,m}^3$ exist. Otherwise only some (or none)
of them may exist.
If any of these cubes does not exists, we let this cube
be the point $\{y\}$.
\end{definition}

Notice that we can only assume that the estimates in \rf{igu} hold for the
cubes $Q$ wich are different from $\{y\}$
(i.e. with $\ell(Q)>0$). So if $\wh{Q}_{y,m}^1 = \{y\}$, say, then, we
only know that $\delta(\wh{Q}_{y,m}^1,2R_0) \leq
m\,A - \alpha_1 -\sigma + \ve_1$.

\begin{lemma}  \label{construc}
Let $y\in \supp(\mu)$. If we choose the constants
$\alpha_1$, $\alpha_2$ and  $A$ big enough, we have
\begin{equation}  \label{inclu}
Q_{y,m}\subset Q_{y,m}^1 \subset \wh{Q}_{y,m}^1 \subset Q_{y,m}^2
\subset \wh{Q}_{y,m}^2\subset Q_{y,m}^3 \subset Q_{y,m-1}.
\end{equation}
\end{lemma}

\begin{proof} Notice  first that for $\alpha_1$, $\alpha_2$ and  $A$ big
enough, then the numbers that appear in the right hand side of the estimates in
\rf{igu} form an estrictly decreasing sequence. That is,
\begin{equation*}
\begin{split}
& m\,A - \ve_1 \,>\, m\,A - \alpha_1 + \ve_1 ,\\
& m\,A - \alpha_1 -\ve_1 \,>\, m\,A - \alpha_1 -\sigma + \ve_1,\\
& m\,A - \alpha_1 -\sigma - \ve_1 \,>\, m\,A - \alpha_1 - \alpha_2 + \ve_1\\
& m\,A - \alpha_1 - \alpha_2 -\ve_1 \,>\, m\,A - \alpha_1 - \alpha_2 -\sigma+\ve_1,\\
& m\,A - \alpha_1 - \alpha_2 -\sigma - \ve_1 \,>\, m\,A - \alpha_1 - \alpha_2
 -2\sigma+\ve_1,\\
& m\,A - \alpha_1 - \alpha_2  -2\sigma-\ve_1 \,>\, (m-1) \,A +  \ve_1.
\end{split}
\end{equation*}

Let us check the inclusion $\wh{Q}_{y,m}^1 \subset Q_{y,m}^2$, for example.
Suppose first that $Q_{y,m}^2\neq \{y\}$, then
$$\delta(Q_{y,m}^2,2R_0) = m\,A - \alpha_1 - \alpha_2 \pm \ve_1.$$
If $\wh{Q}_{y,m}^1 = \{y\}$, the inclusion is obvious. Otherwise,
$$\delta(\wh{Q}_{y,m}^1,2R_0) = m\,A - \alpha_1 -\sigma \pm \ve_1.$$
Then $\delta(\wh{Q}_{y,m}^1,2R_0) >  \delta(Q_{y,m}^2,2R_0)$, and so
$\wh{Q}_{y,m}^1 \subset Q_{y,m}^2$.
Assume now $Q_{y,m}^2 = \{y\}$.
Then,
$$\delta(y,2R_0) \leq m\,A - \alpha_1 - \alpha_2 + \ve_1.$$
In this case there is not any cube $\wh{Q}_{y,m}^1$ satisfying
$$\delta(\wh{Q}_{y,m}^1,2R_0) = m\,A - \alpha_1 -\sigma \pm \ve_1,$$
and so, by our convention, $\wh{Q}_{y,m}^1=\{y\}$. That is, the inclusion
holds in any case.

The other inclusions are proved in a similar way.
\end{proof}

For a fixed $m$, the cubes $Q_{y,m}^1$ may have very diferent sizes for
different $y$'s.
The same happens for the cubes $Q_{y,m}^2$
Nevertheless, in the following lemma we show that we still have some kind
of regularity. This regularity property will be essential for our purposes.

\begin{lemma}  \label{regu}
Let $x,y$ be points in $\supp(\mu)$. Then,
\begin{enumerate}
\item[(a)] If $Q_{x,m}^1\cap Q_{y,m}^1\neq\varnothing$, then
$Q_{x,m}^1\subset \wh{Q}_{y,m}^1$,
in particular $x\in \wh{Q}_{y,m}^1$.
\item[(b)] If $Q_{x,m}^2\cap Q_{y,m}^2\neq\varnothing$,then $Q_{x,m}^2\subset \wh{Q}_{y,m}^2$,
in particular $x\in \wh{Q}_{y,m}^2$.
\end{enumerate}
\end{lemma}

So, although we cannot expect to have the equivalence
$$y\in Q_{x,m}^1 \Leftrightarrow x\in Q_{y,m}^1,$$
we still have something quite close to it, because
the cubes $Q_{x,m}^1$ and $\wh{Q}_{x,m}^1$
are close one each other in the quasimetric $D(\cdot,\cdot)$, since
$\delta(Q_{x,m}^1,\wh{Q}_{x,m}^1)$ is small (at least in front of $A$).
Of course, the same idea applies if we change $1$ by $2$ in the superscripts
of the cubes.

\begin{proof}[Proof of Lemma \ref{regu}]
Let us proof the statement (a). The second statement is proved in an analogous
way. Let $x,y$ be as in (a).
If $\ell(Q_{y,m}^1)>\ell(Q_{x,m}^1)$
(in particular, $Q_{y,m}^1\neq\{y\}$), then $Q_{x,m}^1\subset 3Q_{y,m}^1
\subset \wh{Q}_{y,m}^1$ (the latter inclusion holds provided
$\delta(\wh{Q}_{y,m}^1,2R_0) < \delta(Q_{y,m}^1,2R_0) - 6^nC_0$).

Assume now $\ell(Q_{y,m}^1)\leq \ell(Q_{x,m}^1)$. If $Q_{x,m}^1=\{x\}$,
then $x=y$ and
the result is trivial. If $Q_{x,m}^1\neq\{x\}$, we
denote by $P_{y}$ a cube centered at $y$
with side length $3\ell(Q_{x,m}^1)$. Then, $Q_{x,m}^1\subset P_y\subset
6Q_{x,m}^1$ and so $\delta(Q_{x,m}^1,P_y)\leq 12^n\,C_0$.
Thus
\begin{eqnarray*}
\delta(P_y,2R_0) & \geq & \delta(Q_{x,m}^1,2R_0) - \delta(Q_{x,m}^1,P_y)
- \ve_0\\
& \geq &  \delta(Q_{x,m}^1,2R_0) - 12^n\,C_0 - \ve_0 \\
& \geq &  m\,A - \alpha_1 - \sigma +\ve_1.
\end{eqnarray*}
Therefore,
$\wh{Q}_{y,m}^1\neq\{y\}$ and $\wh{Q}_{y,m}^1\supset P_y\supset
Q_{x,m}^1$.
\end{proof}

Now we are going to define the functions $\vphi_{y,m}$. First we introduce
the auxiliary functions $\psi_{y,m}$.

\begin{definition}  \label{defpsi}
For any $y\in \supp(\mu)\cap 2R_0$,
the function $\psi_{y,m}$ is a function such that
\begin{enumerate}
\item $0\leq \psi_{y,m}(x)\leq \min\left(\dfrac{4}{\ell(Q_{y,m}^1)^n},\,
\dfrac{1}{|y-x|^n}\right)$,
\item $\psi_{y,m}(x) = \dfrac{1}{|x-y|^n}$ if $x\in \wh{Q}_{y,m}^2\setminus
Q_{y,m}^1$,
\item $\supp(\psi_{y,m})\subset Q_{y,m}^3$,
\item $|\psi_{y,m}'(x)|\leq C_{12}\,
\min\left(\dfrac{1}{\ell(Q_{y,m}^1)^{n+1}},\,  \dfrac{1}{|y-x|^{n+1}}\right)$.
\end{enumerate}
\end{definition}

It is not difficult to check that such a function exists if we  choose $C_{12}$ big
enough.
We have to take into account that $2\wh{Q}_{y,m}^2 \subset Q_{y,m}^3$.
This is due to the fact that $\delta(\wh{Q}_{y,m}^2,2\wh{Q}_{y,m}^2)\leq 4^nC_0
< \delta(\wh{Q}_{y,m}^2,Q_{y,m}^3)$ if $\ell(\wh{Q}_{y,m}^2)\neq0$.

In the definition of $\psi_{y,m}$, if $Q_{y,m}^1=\{y\}$, then one must take
$1/\ell(Q_{y,m}^1)=\infty$. If $\wh{Q}_{y,m}^2 = \{y\}$, then we set
$\psi_{y,m}\equiv0$. This choice satisfies the conditions for the definition
of $\psi_{y,m}$ stated above.

Choosing $\alpha_2$ big enough, the largest part of the $L^1(\mu)$ norm of
$\psi_{y,m}$ will come from the integral over $Q_{y,m}^2\setminus
\wh{Q}_{y,m}^1$.
We state this in a precise way in the following lemma.

\begin{lemma}  \label{propor}
There exists some constant $\ve_2$ depending on $n$, $d$, $C_0$, $\ve_0$,
$\ve_1$ and $\sigma$ (but not on $\alpha_1$, $\alpha_2$ nor $A$) such that
if $Q_{y,m}^1 \neq \{y\}$, then
\begin{equation}  \label{err30}
\left| \|\psi_{y,m}\|_{L^1(\mu)} - \alpha_2 \right| \leq \ve_2
\end{equation}
and
\begin{equation}  \label{err31}
\left| \|\psi_{y,m}\|_{L^1(\mu)} - \int_{Q_{y,m}^2\setminus \wh{Q}_{y,m}^1}
 \frac{1}{|y-x|^n}\,d\mu(x) \right| \leq \ve_2.
\end{equation}
\end{lemma}

The proof of this result is an easy calculation that we will skip.
A direct consequence of it is
$$\lim_{\alpha_2\to \infty} \frac{1}{\alpha_2}\int_{Q_{y,m}^2\setminus
\wh{Q}_{y,m}^1}
\frac{1}{|y-x|^n}\,d\mu(x) = 1$$
for $y\in\supp(\mu)$ such that $\delta(y,2R_0)>m\,A$.

\begin{definition}  \label{defphi}
Let $w_{i,m}$ be the weight function defined for $y\in \bigcup_{i\in I^\prime_m}
Q_{i,m}$ (these are the cubes of $\DD_m$ with $\ell(Q_{i,m})>0$) by
$$w_{i,m}(y) = \frac{\chi_{Q_{i,m}}(y)}{\sum_{j\in I^\prime_m}
\chi_{Q_{j,m}(y)}}.$$
If $y\in \supp(\mu)\cap 2R_0$ belongs to some cube $Q_{i,m}$ centered at some
point $y_i$, with
$\ell(Q_{i,m})>0$, then we set
$$\vphi_{y,m}(x) = \alpha_2^{-1}\, \sum_i w_{i,m}(y)\,\psi_{y_i,m}(x).$$
If $y$ does not belong to any cube $Q_{i,m}$ with
$\ell(Q_{i,m})>0$ (this implies $\delta(y,2R_0)\leq mA$
and $Q_{y,m} = \{y\}$), then we set
$$\vphi_{y,m}(x) = \alpha_2^{-1}\, \psi_{y,m}(x).$$
\end{definition}

Setting $w_{i,m}(y) = \chi_{Q_{i,m}}(y)$ if $\ell(Q_{i,m})=0$, we can write
$$\vphi_{y,m}(x) = \alpha_2^{-1}\, \sum_i w_{i,m}(y)\,\psi_{y_i,m}(x),$$
for {\em any} $y$.

Let us remark that a more natural definition for $\vphi_{y,m}$
would have been the choice
$\vphi_{y,m}(x) = \alpha_2^{-1}\, \psi_{y,m}(x)$ for all $y$. However, as we
shall see, for some of the arguments in the proof of the Main Lemma below
(in Subsection \ref{correccio}), the
choice of Definition \ref{defphi} is better.

In order to study some of the properties of the functions
$\vphi_{y,m}$, we need to introduce some additional notation.

\begin{definition}
Given $x\in\supp(\mu)$, we denote by $\QH_{x,m}$ a doubling cube centered at
$x$ such that $\delta(\QH_{x,m},2R_0)= m\,A - \alpha_1- \alpha_2 - 3\,\sigma
\pm \ve_1$.
Also, we denote by $\check{Q}_{x,m}^1$ and $\Check{\Check{Q}}_{x,m}^1$ some
doubling cubes centered at $x$ such that
\begin{equation*}
\begin{split}
& \delta(\check{Q}_{x,m}^1,2R_0) = m\,A - \alpha_1 + \sigma \pm \ve_1,\\
& \delta(\Check{\Check{Q}}_{x,m}^1,2R_0) = m\,A - \alpha_1 + 2\sigma \pm \ve_1
\end{split}
\end{equation*}
(the idea is that the symbols $\,\,\wh{ }\,\,$ and $ \,\,\check{ }\,\,$
are inverse operations, modulo some small errors).
If any of the cubes
$\check{Q}_{x,m}^1,\Check{\Check{Q}}_{x,m}^1, \QH_{x,m}$
does not exist, then we let it be the point $x$.
\end{definition}

So, when $\delta(x,2R_0)$ is big enough, one should think that $\QH_{x,m}$ is
a cube a little bigger than
$\wh{Q}^3_{x,m}$, while $\check{Q}_{x,m}^1$ is a little smaller than
$Q_{x,m}^1$. Also, $\Check{\Check{Q}}_{x,m}^1$ is a little smaller than
$\check{Q}_{x,m}^1$, but still much bigger than $Q_{x,m}$.

\begin{lemma} \label{phi}
Let $x,y\in \supp(\mu)$. For $\alpha_1$ and $\alpha_2$ big enough, we have:
\begin{itemize}
\item[(a)] If $x\in Q_{x_0,m}$ and $y\not\in \wh{Q}_{x_0,m}^3$, then
$\vphi_{y,m}(x) = 0$. In particular,
$\vphi_{y,m}(x) = 0$ if $y\not\in \wh{Q}_{x,m}^3$.

\item[(b)] If $y\in \check{Q}_{x,m}^1$, then $\ds \vphi_{y,m}(x) \leq C\,
\frac{\alpha_2^{-1}}{\ell(\check{Q}_{x,m}^1)^n}.$

\item[(c)] Let $\ve_3>0$ be an arbitrary constant. If
$\alpha_1$ is big enough (depending on $\ve_3,\,C_0,\,n,\,d$ but not on 
$\alpha_2$), then
$$\vphi_{y,m}(x) \leq \frac{\alpha_2^{-1}\, (1+\ve_3/2)}{|y-x|^n}\qquad
\mbox{if $y\not\in \check{Q}_{x,m}^1$,}$$
and
$$\ds \vphi_{y,m}(x) \geq \frac{\alpha_2^{-1}\,(1-\ve_3/2)}{|y-x|^{n}}
\qquad \mbox{if $y\in Q_{x,m}^2 \setminus \wh{Q}_{x,m}^1$.}$$

\item[(d)] If $x\in Q_{x_0,m}$, then $$|\ds \vphi_{y,m}'(x)| \leq
C\, \alpha_2^{-1} \,
\min\left(\frac{1}{\ell(\check{Q}_{x_0,m}^1)^{n+1}}, \, \frac{1}{|y-x|^{n+1}}
\right).$$
\end{itemize}
\end{lemma}

Notice that, in Definition \ref{defpsi} of the functions $\psi_{y,m}$,
the properties that define these functions are stated with respect to cubes
centered at $y$ ($Q_{y,m}^1$, $Q_{y,m}^2$, $Q_{y,m}^3$...). In this lemma
some analogous properties are stated, but these properties have to do with
cubes centered at $x$ or containing $x$
($Q_{x_0,m}$, $\check{Q}_{x,m}^1$, $Q_{x,m}^2$, $\wh{Q}_{x,m}^3$...).

\begin{proof}
\begin{itemize}

\item[(a)]
Let $x_0\in \supp(\mu)$ and $x\in Q_{x_0,m}$.
If $\vphi_{y,m}(x)\neq0$, there exists some $i$ with $y\in Q_{i,m} \equiv
Q_{y_i,m}$ and
$x\in Q_{y_i,m}^3$. Then $Q_{x_0,m}^3 \cap Q_{y_i,m}^3\neq\varnothing$
and so
$y\in Q_{y_i,m}^3\subset \wh{Q}_{x_0,m}^3$ (as in Lemma \ref{regu}).

\vv
\item[(b)] Let $y\in \check{Q}_{x,m}^1$ and
let $y_i$ be such that $y\in Q_{y_i,m}$. We know that
$$\vphi_{y_i,m}(x) \leq C\,\alpha_2^{-1}\,
\frac{1}{\ell(Q_{y_i,m}^1)^{n}}.$$
So we are done if we see that $\ell(Q_{y_i,m}^1) \geq \ell(\check{Q}_{x,m}^1)$.

As in Lemma \ref{regu}, we have
$$y\in \check{Q}_{x,m}^1 \Rightarrow \check{Q}_{y_i,m}^1 \cap\check{Q}_{x,m}^1
\neq \varnothing
\Rightarrow \check{Q}_{x,m}^1\subset Q_{y_i,m}^1.$$
Thus $\ell(\check{Q}_{x,m}^1)\leq \ell(Q_{y_i,m}^1).$

\vv
\item[(c)]
Let us see the first inequality.
If $y\not\in\check{Q}_{x,m}^1$ and $y$ belongs to some cube $Q_{y_i,m}$
with $\ell(Q_{y_i,m})>0$,
then $x\not\in\Check{\Check{Q}}_{y_i,m}^1$ because
otherwise, as in Lemma \ref{regu}, we would get
$\Check{\Check{Q}}_{y_i,m}^1 \subset\check{Q}_{x,m}^1.$
However, since we assume $\alpha_1\gg \sigma$, the cube
$\Check{\Check{Q}}_{y_i,m}^1$ is bigger than $Q_{y_i,m}$ and contains $y$.
So $y\in\check{Q}_{x,m}^1$, which is a contradiction.

Since $x\not\in\Check{\Check{Q}}_{y_i,m}^1$ and this cube is much bigger than
$Q_{y_i,m}$, if $\alpha_1$ is big enough we get
$$\frac{\alpha_2^{-1}}{|y_i-x|^n} \leq
\frac{\alpha_2^{-1}\,(1+\ve_3)}{|y-x|^n}.$$
As this holds for all $i$ with $w_{i,m}(y)\neq0$, we obtain
$$\vphi_{y,m}(x) \leq \frac{\alpha_2^{-1}\,(1+\ve_3)}{|y-x|^n}.$$
This inequality also holds if $\ell(Q_{y_i,m})=0$ with $\ve_3=0$, since in this
case $y_i=y$.

\vv
We consider now the second inequality in (c).
Let $y\in\supp(\mu)$ be such that
$y\in Q_{x,m}^2\setminus \wh{Q}_{x,m}^1.$
If $y\in Q_{y_i,m}$ with $\ell(Q_{y_i,m})>0$ for some $i$,
by Lemma \ref{regu} we get
$x\in \wh{Q}_{y_i,m}^2\setminus Q_{y_i,m}^1.$
Since this is satisfied for all $i$ such that $w_{i,m}(y)\neq 0$,
$$\vphi_{y,m}(x) = \sum_i w_{i,m}(y) \frac{\alpha_2^{-1}}{|y_i-x|^n}.$$
If $\alpha_1$ has been chosen big enough, then $\ell(Q_{y_i,m}^1) \gg \ell(Q_{y_i,m})$
and one has
$$\frac{\alpha_2^{-1}}{|y_i-x|^n} \geq
\frac{\alpha_2^{-1}\,(1-\ve_3/2)}{|y-x|^n}.$$
Thus
\begin{equation}  \label{uiop}
\vphi_{y,m}(x) \geq \frac{\alpha_2^{-1}\,(1-\ve_3/2)}{|y-x|^n}.
\end{equation}
If $y\in Q_{x,m}^2\setminus \wh{Q}_{x,m}^1$ and
$y\in Q_{i,m}$ with $\ell(Q_{i,m})=0$, then by Lemma \ref{regu} we also get
$x\in \wh{Q}_{y,m}^2\setminus Q_{y,m}^1$
(in particular $\wh{Q}_{y,m}^2\neq \{y\}$). Then \rf{uiop} holds in this
case too (with $\ve_3=0$).

\vv
\item[(d)]
Suppose first that $y\in \check{Q}_{x_0,m}^1$. In this case we must show that
$$|\vphi_{y,m}'(x)| \leq
C\, \frac{\alpha_2^{-1}}{\ell(\check{Q}_{x_0,m}^1)^{n+1}}.$$
Let $y_i$ be such that $y\in Q_{y_i,m}$. We know that
$$|\vphi_{y_i,m}'(x)| \leq C\,
\frac{\alpha_2^{-1}}{\ell(Q_{y_i,m}^1)^{n+1}}.$$
By the definition of $\vphi_y(x)$, it is enough to see that
$\ell(Q_{y_i,m}^1) \geq \ell(\check{Q}_{x_0,m}^1)$. This follows from the inclusion
$Q_{y_i,m}^1 \supset \check{Q}_{x_0,m}^1$, which holds because
$y\in \check{Q}_{y_i,m}^1 \cap \check{Q}_{x_0,m}^1$ and then we can apply Lemma
\ref{regu} (in fact, a slight variant of Lemma \ref{regu}).

Suppose now that $y\not\in \check{Q}_{x_0,m}^1$. It is enough to show that
$$|\vphi_{y,m}'(x)| \leq
C\, \frac{ \alpha_2^{-1}}{|y-x|^{n+1}}.$$
Let $y_i$ be such that $y\in Q_{y_i,m}$. By definition we have
$$|\vphi_{y_i,m}'(x)| \leq C\,\frac{\alpha_2^{-1}}{|y_i-x|^{n+1}}.$$
We are going to see that
\begin{equation} \label{rar0}
|y-y_i|\leq |y-x|/2.
\end{equation}
Assume $|y-y_i| > |y-x|/2$. Then, since $x\in \frac{1}{2}\check{Q}_{x_0,m}^1$
(for $\alpha_1$ big enough),
\begin{equation} \label{rar}
\ell(Q_{y_i,m})> C^{-1}\,|y-x| \geq C^{-1}\,
\ell(\check{Q}_{x_0,m}^1).
\end{equation}
Notice that from the first inequality in \rf{rar} we get $\dist(x,Q_{y_i,m})
\leq
C\,\ell(Q_{y_i,m})$. In this situation we have
$\check{Q}_{x_0,m}^1\subset C\,Q_{y_i,m} \subset
\Check{\Check{Q}}_{y_i,m}^1.$
This is not possible, since by Lemma \ref{regu} we would have
$\check{Q}_{x_0,m}^1\supset \Check{\Check{Q}}_{y_i,m}^1,$
and then we would get $\check{Q}_{x_0,m}^1 = \Check{\Check{Q}}_{y_i,m}^1$.
This would imply $x_0= y_i$ and also
$x_0= y_i=\check{Q}_{x_0,m}^1 = \Check{\Check{Q}}_{y_i,m}^1$, and then $y=y_i$
which is a contradiction because we are assuming that \rf{rar0} does not hold.

So \rf{rar0} is true and $|y_i-x|\approx |y-x|$. Thus
$$|\vphi_{y_i,m}'(x)| \leq C\,
\frac{\alpha_2^{-1}}{|y-x|^{n+1}}.$$
Since this holds for any $i$ such that $y\in Q_{y_i,m}$, we get
$$|\vphi_{y,m}'(x)| \leq C\, \frac{\alpha_2^{-1}}{|y-x|^{n+1}}.$$

\end{itemize}
\end{proof}

Some of the estimates in the preceding lemma will be used to prove
next result, which was one of our main goals in this section.

\begin{lemma} \label{convo}
For any $\ve_3>0$, if $\alpha_1$ and $\alpha_2$ are big
enough, for all $x\in \supp(\mu)$ we have
\begin{equation} \label{convo1}
\int \vphi_{y,m}(x)\,d\mu(y) \leq 1 + \ve_3.
\end{equation}
If $x\in\supp(\mu)$ is such that there exists some cube $Q\in\DD_m$ with $Q\ni
x$ and $\ell(Q)>0$ (in particular if $\delta(x,2R_0)>m\,A$), then
\begin{equation} \label{convo2}
1-\ve_3 \leq \int \vphi_{y,m}(x)\,d\mu(y)
\end{equation}
\end{lemma}

Let us observe that if $\mu$ were invariant by translations and
$\vphi_{y,m}(x)=\vphi_{m}(y-x)$, then \rf{convo1} and  \rf{convo2}
would hold with $\ve_3=0$ (choosing $\|\vphi_{y,m}\|_{L^1(\mu)} $ $=1$).

\begin{proof}
Let us see  \rf{convo2} first. So we assume
that there exist some cube $Q_{i,m}\in \DD_m$ containing $x$ with
$\ell(Q_{i,m})>0$. Since $x\in Q_{i,m}\subset \check{Q}_{i,m}^1$, we have
$\check{Q}_{i,m}^1 \subset Q_{x,m}^1$. In particular, $\ell(Q_{x,m}^1)>0$.
By Lemma \ref{propor} and the second inequality of (c) in Lemma \ref{phi}
we get
\begin{eqnarray*}
\int \vphi_{y,m}(x)\,d\mu(y) & \geq & \int_{Q_{x,m}^2\setminus \wh{Q}_{x,m}^1}
 \vphi_{y,m}(x)\,d\mu(y)\\
& \geq & \int_{Q_{x,m}^2\setminus \wh{Q}_{x,m}^1}\frac{\alpha_2^{-1}\,(1-\ve_3/2)}{|y-x|^n}
\,d\mu(y)\\
& \geq & \alpha_2^{-1}\,(\alpha_2-2\ve_2)\,(1-\ve_3/2).
\end{eqnarray*}
So \rf{convo2} holds if we take $\alpha_2$ big enough.

Consider now \rf{convo1}.
By (a) in Lemma \ref{phi} have
$$\int \vphi_{y,m}(x)\,d\mu(y) = \int_{\wh{Q}_{x,m}^3}
\vphi_{y,m}(x)\,d\mu(y).$$
Thus we can write
\begin{equation}  \label{desc}
\int \vphi_{y,m}(x)\,d\mu(y) = \int_{\wh{Q}_{x,m}^3\setminus\check{Q}_{x,m}^1}
\vphi_{y,m}(x)\,d\mu(y) + \int_{\check{Q}_{x,m}^1} \vphi_{y,m}(x)\,d\mu(y).
\end{equation}

Let us estimate the first integral on the right hand side of \rf{desc}.
Using the first inequality in (c) of Lemma \ref{phi} we obtain
\begin{eqnarray}  \label{desc3}
\int_{\wh{Q}_{x,m}^3\setminus\check{Q}_{x,m}^1} \vphi_{y,m}(x)\,d\mu(y)
& \leq & \int_{\wh{Q}_{x,m}^3\setminus\check{Q}_{x,m}^1}
\frac{\alpha_2^{-1}\,(1+\ve_3/2)}{|y-x|^n}\, d\mu(y) \nonumber \\
& = & \delta(\check{Q}_{x,m}^1,\wh{Q}_{x,m}^3) \,\alpha_2^{-1}\,(1+\ve_3/2)
\nonumber\\
& \leq & \alpha_2^{-1}\,
(\alpha_2 + 4\,\sigma + 2\,\ve_1)\,(1+\ve_3/2).
\end{eqnarray}

Let us consider the last integral in \rf{desc}
(only in the case $\check{Q}_{x,m}^1\neq\{x\}$).
By (b) in Lemma \ref{phi} we have
\begin{equation} \label{desc4}
\int_{\check{Q}_{x,m}^1} \vphi_{y,m}(x)\,d\mu(y)
\leq \int_{\check{Q}_{x,m}^1}
\frac{C\,\alpha_2^{-1}}{\ell(\check{Q}_{x,m}^1)^n}\,d\mu(y)
\leq C\,C_0\,\alpha_2^{-1}.
\end{equation}
From \rf{desc3} and \rf{desc4} we get \rf{convo1}.
\end{proof}


\section{Proof of the Main Lemma}  \label{seccore}


\subsection{The argument}

As stated above, $A$ is a large positive integer that will be fixed at the end
of the proof.
We assume that the support of $f$ is contained in some doubling cube $R_0$, and
for each integer $m\geq1$ we consider the family $\DD_m$ of ``dyadic''
cubes $Q_{i,m}$, $i\in I_m$,
introduced in Definition \ref{defcubs}, and we set $\DD =
\bigcup_{m\geq1} \DD_m$. Recall that the elements of $\DD$ may be cubes with
side length $0$, i.e. points.

For each $m$ we will construct functions $g_m$ and $b_m$.
The function $g_m$ will be supported on a subfamily $\DD_m^G$ of the cubes in
$\DD_m$. On the other hand, $b_m$
will be supported on a subfamily $\DD_m^B$ of the cubes in $\DD_m$.
We set $\DD^G =
\bigcup_{m\geq1} \DD_m^G$ and $\DD^B = \bigcup_{m\geq1}\DD_m^B$.
The cubes in $\DD^G$ will be called good cubes and the ones in $\DD^B$ bad
cubes (let us remark that in the family $\DD_m$, in general, there are also
cubes which are neither good nor bad).

From $g_m$ and $b_m$, we will obtain the following potentials:
\begin{eqnarray*}
U_m^G(x) & = & \int \vphi_{y,m}(x)\,g_m(y)\,d\mu(y),\\
U_m^B(x) & = & \int \vphi_{y,m}(x)\,b_m(y)\,d\mu(y), \\
U_m(x) & = & U_m^G(x) + U_m^B(x).
\end{eqnarray*}
This potentials will be successively subtracted from $f$. We will set
$$f_{m+1}(x) = f(x) - \sum_{j=1}^{m} U_j(x) = f_m(x) - U_m(x)$$
and
\begin{equation}  \label{hh00}
h_0 = f - \sum_{m=1}^\infty U_m = \lim_{m\to\infty} f_m.
\end{equation}
The support of the functions $g_m$, $b_m$, $U^G_m$, $U^B_m$ will be contained
in $2R_0$.

By induction we will show that
the functions $g_m$, $b_m$, $U_m$ and $f_m$ fulfil the following
properties:


\begin{itemize}
\item[(a)] $|g_m|,\, |b_m|\leq C_8\, A\,\|f\|_*$.

\item[(b)] $|m_Q f_{m+1}| \leq A \,\|f\|_*$ if $Q\in \DD_{m}$ and $\ell(Q)>0$.

\item[(c)] If $g_m\not\equiv 0$ on $Q$, $Q\in \DD_m$, with $\ell(Q)>0$,
then $|m_Q f_{m+1}| \leq \dfrac{7}{20}\,A\,\|f\|_*$.

\item[(d)] If $Q\in \DD_m$ and $|m_Q f_m|\leq \dfrac{8}{20}\,A\,\|f\|_*$, then
$U_m\equiv 0$ and $g_m\equiv b_m\equiv 0$ on $Q$.

\item[(e)] If $Q\in \DD_m$ and $\delta(Q,2R_0)\leq (m-\frac{1}{10})\,A$
(so $\ell(Q)=0$), then
$U_m\equiv 0$ and $g_m\equiv b_m\equiv 0$ on $Q$.

\end{itemize}

\vv
Finally, we will see that our construction satisfies the following properties too:

\begin{itemize}

\item[(f)] If $\delta(x,2R_0)<\infty$, then $|h_0(x)|\leq C_{9}\,A\,\|f\|_*,$
and if $Q\in\DD_m$ and $\ell(Q)=0$, then $|m_Q f_{m+1}| \equiv |f_{m+1}(z_Q)|
\leq C_{9}\,A\,\|f\|_*$.

\item[(g)] For each $m$, there are functions $g^1_m,\ldots,g_m^N$ such
that
\begin{itemize}
\item[(g.1)]
$\ds U_m^G(x) = \sum_{p=1}^N\int \vphi_{y,m}^p(x)\,g_m^p(y)\,d\mu(y),$
where $\vphi_{y,m}^p$ is defined below.

\item[(g.2)]
$|g_m^p|\leq 2C_8\, A\,\|f\|_*$ for $p=1,\ldots,N,$

\item[(g.3)] The functions $\sum_{p=1}^N |g_m^p|$ have disjoint supports for
different $m$'s.
\end{itemize}

\item[(h)] The family of cubes $\DD^B$ that support the functions $b_m$,
$m\geq1$, satisfies the following Carleson packing condition for each
cube $R\in \DD_m$ with $\ell(R)>0$:
\begin{equation}  \label{pack}
\sum_{\substack{ Q:\,Q\cap R\neq\varnothing \\Q\in D^B_k,\,k> m}}
\mu(Q) \leq C\,\mu(R).
\end{equation}

\end{itemize}


Let us remark that if some cube $Q$ coincides with a point $\{x\}$, then we set
$m_Q f_m\equiv f_m(x)$. Also, the notation for the sum in (h) is an abuse
of notation. This sum has to be undestood as
$$\sum_{\substack{ Q:\,Q\subset 2R\\Q\in D^B_k,\,k>m}}
\mu(Q) \equiv \sum_{\substack{ Q:\,\ell(Q)>0,\,Q\subset 2R\\Q\in D^B_k,\,k>m}}
\mu(Q) + \sum_{k > m} \mu\left\{x\in 2R:\,\{x\}\in D_k^B\right\}.$$
On the other hand, the number $N$ that appears in (g) is the
number of disjoint
families of cubes given in the  Covering  Theorem of Besicovich,
which only depends only on $d$.

The functions $\vphi_{y,m}^p$ of (g) are defined as follows. 
We set $\DD_m = \DD_m^{1} \cup \cdots \cup\DD_m^{N}$,
where each subfamily $\DD_m^{p}$ is disjoint (recall that the cubes of $\DD_m$
originated from Besicovich's Covering Theorem). Then we set
$$\vphi_{y,m}^p(x) = \vphi_{y_i,m}(x)$$
if $y\in Q_{i,m}$ with $Q_{i,m}\in \DD_m^{p}$, and $\vphi_{y,m}^p(x)\equiv0$ if there
does not exist any cube of the subfamily $\DD^p_m$ containing $y$.

\vspace{2mm}
First we will show that if there exist functions $g_m$ and $b_m$ satisfying
(a)--(h) then the Main Lemma follows, and later we
will show the existence of these functions.


It is not difficult to check that if \rf{main} and \rf{main2} hold,mi
then the sum of \rf{hh00} converges in $L^1_{loc}(\mu)$ (this is left
to the reader). Since the support of all the functions involved is contained in
$2R_0$, the convergence is in $L^1(\mu)$.

Let us see now that if (b) and (f) hold, then $\|h_0\|_{L^\infty(\mu)}
\leq C\,A\, \|f\|_*$.
Taking into account (f), we only have to see that $|h_0(x)| \leq C\,A\, \|f\|_*$
for $x\in\supp(\mu)$ such that $\delta(x,2R_0)=\infty$.
In this case, if $Q\in\DD_k$ is such that $x\in Q$, then $\ell(Q)>0$.
We are going to see that
\begin{equation}  \label{uuj}
|m_Q f_m| \leq C\,A \,\|f\|_*\qquad\mbox{for $Q\in \DD_{k}$, $k\leq m-1$}
\end{equation}
(not only for $k=m-1$, which is a direct consequence of (b) and (f)).
Take $Q\in D_k$, $k< m-1$. This cube is covered with finite overlap
by the family of cubes $\DD_{m-1}$. Moreover, if $P\in \DD_{m-1}$ and $P\cap
Q\neq\varnothing$, then $\ell(P)\leq \ell(Q)/10$ by Lemma \ref{mides}, and so
$P\subset 2Q$. Thus we get
$$\int_Q |f_m|\,d\mu \leq \sum_i \int_{Q\cap Q_{i,m-1}} |f_m| \,d\mu
\leq C\,A\,\|f\|_*\,\mu(2Q) \leq C\,A\,\|f\|_*\,\mu(Q),$$
and \rf{uuj} follows (notice that, as remarked above, we have abused notation
for the cubes which are single points).

Then $h_0$ will satisfy $|m_Q h_0| \leq C\,A \,\|f\|_*$ for all $Q\in \DD$
containing $x$,
because the sequence $\{f_m\}_m$ converges to $h_0$ in $L^1(\mu)$. Then, by the
Lebesgue differentiation theorem we will get that $|h_0(x)| \leq
C\,A \,\|f\|_*$ (this theorem can be applied to the cubes $Q\in\DD$ which
are non centered because they are doubling) for $\mu$-a.e. $x\in\supp(\mu)$
with $\delta(x,2R_0)=\infty$. Therefore, $\|h_0\|_{L^\infty(\mu)} \leq
C\,A \,\|f\|_*$.


\vv
Observe that the functions $g_m^p$ in (g.1) originate the same potential
as $g_m$. In fact, they will be constructed modifying slightly the function
$g_m$ in such a way that they are supported in disjoint sets for different
$m$'s. By (g.2)  we have
$$\sum_m \sum_{p=1}^N |g_m^p| \leq 2N\,C_8\,A\,\|f\|_*.$$

The supports of the functions $b_m$ may be not disjoint.
To solve this problem, we will construct ``corrected'' versions ($b_m^p$,
$p=1,\ldots,N$) of $w_{i,m}\, b_m$.
Moreover, as in the case of $g_m$, the modifications will be
made in such a way that the potentials $U_m^B$
will not change.


\subsection{The ``correction'' of $b_m$} \label{correccio}

We assume that the functions $b_m$, $m\geq1$, have been obtained and they
satisfy (a)--(h). We will start the construction of some new functions
(the corrected versions of $w_{i,m}\,b_m$)
in the small cubes, and then we will go over
the cubes from previous generations. However, since there is an infinite
number of generations, we will need to use a limiting argument.

For each $j$ we can write the potential originated by $b_j$ as
$$U^B_j(x) = \sum_{i\in I_j} \vphi_{y_i,j}(x) \int w_{i,j}(y)\, b_j(y)\,d\mu(y).$$
For a fixed $m\geq 1$ we are going to define functions $v_{i,j}^m$, for
$j=m$, $m$\nobreakdash$-1, \,\ldots, 1$, and all $i\in I_j$. The functions
$v_{i,j}^m$
will satisfy
\begin{equation}  \label{cor0}
\supp(v_{i,j}^m) \subset Q_{i,j},
\end{equation}
where $Q_{i,j}\in \DD_J^B$,
the sign of $v_{i,j}^m$ will be constant on $Q_{i,j}$,
and
\begin{equation}  \label{cor1}
\int v_{i,j}^m(y)\,d\mu(y) = \int w_{i,j}(y)\, b_j(y)\,d\mu(y).
\end{equation}
Moreover, we will also have
\begin{equation}  \label{cor2}
\sum_{j=1}^m \sum_{i\in I_j} |v_{i,j}^m|\leq C_{11}\,A\,\|f\|_*.
\end{equation}

We set $v_{i,m}^m(y)= w_{i,m}(y)\, b_m(y)$ for all $i\in I_m$.
Assume that we have obtained functions $v_{i,m}^m,\, v_{i,m-1}^m,\ldots
v_{i,k+1}^m$ for all the $i$'s, fulfiling \rf{cor0}, \rf{cor1},
and such that
$$\sum_{j=k+1}^m \sum_{i\in I_j}|v_{i,j}^m| \leq B\,A\,\|f\|_*,$$
where $B$ is some constant that will be fixed below.
We are going to construct $v_{i,k}^m$ now.

Let $Q_{i_0,k}\in\DD_k$ be some fixed cube from the $k$-th generation. Assume
first that $Q_{i_0,k}$ is not a single point.
Since the cubes in the family $\DD^B$ satisfy the
packing condition \rf{pack},
for any $t>0$ we get
\begin{eqnarray*}
\lefteqn{
\mu\biggl\{y\in Q_{i_0,k}:\, \sum_{j=k+1}^m \sum_{i\in I_j} |v_{i,j}^m(y)| >t
\biggr\} } &&\\
& \leq &
\frac{1}{t}  \sum_{j=k+1}^m \sum_{i\in I_j} \int_{Q_{i_0,k}} |v_{i,j}^m(y)|\,
d\mu(y)
\\ & \leq &
\frac{1}{t}\sum_{j=k+1}^m \sum_{i\in I_j}
\int_{Q_{i_0,k}} |w_{i,j}(y)\, b_j(y)|\,d\mu(y) \\
& \leq &
\frac{C_8\,A\,\|f\|_*}{t}
\sum_{\substack{ Q:\,Q\cap Q_{i_0,k}\neq \varnothing
\\Q\in D^B_j,\,j> k}} \mu(Q)
\leq \frac{C_{12}\,A\,\|f\|_*}{t}\,\mu(Q_{i_0,k}).
\end{eqnarray*}
Therefore, if we choose $t= 2C_{12}\,A\,\|f\|_*$ and we denote
$$V_{i_0,k}^m =\biggl\{y\in Q_{i_0,k}:\, \sum_{j=k+1}^m \sum_{i\in I_j}
|v_{i,j}^m(y)| \leq t \biggr\},$$
we have
$\mu(V_{i_0,k}^m) \geq \frac{1}{2}\mu(Q_{i_0,k})$.
If we set
$v_{i_0,k}^m = c_{i_0,k}^m\,\chi_{V_{i_0,k}},$
where $c_{i_0,k}^m\in\R$ is such that \rf{cor1} holds for $i=i_0$,
then
$$|c_{i_0,k}^m| \leq \frac{1}{\mu(V_{i_0,k})}
\int |w_{i,k}(y)\, b_k(y)|\,d\mu(y) \leq 2\,C_8\,A\,\|f\|_*.$$
By the finite overlap of the cubes in $\DD_k^B$, we get
$$
\sum_{\substack{i_0:\, Q_{i_0,k}\in \DD_k^B\\  \ell(Q_{i_0,k})\neq 0}}
|v_{i_0,k}^m| \leq C_B\,2\,C_8\,A\,\|f\|_*,
$$
where $C_B$ is the overlap constant in the Covering Theorem of Besicovich.
Now if we take $B:=2\,C_B\,C_8 + 2C_{12}$, we will have
\begin{equation}  \label{overl}
\sum_{\substack{i_0:\, Q_{i_0,k}\in \DD_k^B\\  \ell(Q_{i_0,k})\neq 0}}
|v_{i_0,k}^m| +
\sum_{j=k+1}^m \sum_{i\in I_j}|v_{i,j}^m| \leq B\,A\,\|f\|_*.
\end{equation}

In case $Q_{i_0,k}$ is a single point $\{y\}$, then we set
$v_{i_0,k}^m(y)= w_{i_0,k}(y)\,b_k(y) = b_k(y)$. All the cubes of
the generations $k+1,\ldots, m$ that intersect $Q_{i_0,k}\equiv \{y\}$ coincide
with $\{y\}$ by Lemma \ref{mides}. From (e) we get that
$b_{k+1}(y)=b_{k+2}(y)= \cdots =0$, which is the same as saying
that $v_{i,k+1}^m(y)=v_{i,k+2}^m(y)=\cdots=0$ for all $i$. So we have
\begin{equation}  \label{overll}
\sum_{j=k}^m \sum_{i\in I_j}|v_{i,j}^m(y)| = |b_k(y)| \leq C_8\,A\,\|f\|_*
\leq B\,A\,\|f\|_*.
\end{equation}
From \rf{overl} and \rf{overll} we get
$$\sum_{j=k}^m \sum_{i\in I_j}|v_{i,j}^m| \leq B\,A\,\|f\|_*.$$
Operating in this way, the functions $v_{i,j}^m$,
$j=m, \,m-1,\ldots, 1$, $i\in I_j$, will satisfy the conditions
\rf{cor0}, \rf{cor1} and \rf{cor2} (with $C_{11}=B$).

Now we can take a subsequence $\{m_k\}_k$ such that for all $i\in I_1$ (i.e.
for all the cubes of the first generation) the functions $\{v_{i,1}^{m_k}\}_k$
converge weakly in $L^\infty(\mu)$ to some function $v_{i,1}\in L^\infty(\mu)$.
Let us remark that the sequence $\{m_k\}_k$ can be chosen independently of $i$
since, by the Besicovich's Covering Theorem, there is a bounded number $N$
of subfamilies $\DD_1^{1}, \ldots,\DD_1^{N}$
of $\DD_1$ such that each subfamily $\DD_1^{p}$ is disjoint. If we denote by
$\DD_1^{p,B}$ the subfamily of bad cubes of $\DD_1^p$, we can write
$$\sum_{i\in I_1} v_{i,1}^m = \sum_{p=1}^N \sum_{i:\,Q_{i,1}\in\DD_1^{p,B}}
v_{i,1}^m,$$
and we can choose $\{m_k\}_k$ such that, for each $p$,
$\sum_{i:\,Q_{i,1}\in\DD_1^{p,B}} v_{i,1}^{m_k}$ converges weakly to
$\sum_{i:\,Q_{i,1}\in\DD_1^{p,B}} v_{i,1}$.

In a similar way, we can consider another subsequence of $\{m_{k_j}\}_j$ of
$\{m_k\}_k$ such that for all $i\in I_2$ the functions $\{v_{i,2}^{m_{k_j}}\}_j$
converge weakly in $L^\infty(\mu)$ to some function $v_{i,2}\in L^\infty(\mu)$.
Going on with this process, we will obtain functions $v_{i,j}$, $j\geq1$, that
satisfy \rf{cor0}, \rf{cor1} (without the
superscript $m$) and
\begin{equation}  \label{cor2'}
\sum_{j=1}^\infty\sum_{i\in I_j} |v_{i,j}|\leq C_{11}\,A\,\|f\|_*.
\end{equation}
Also, we have
$$U^B_j(x) = \sum_{i\in I_j} \vphi_{y_i,j}(x) \int v_{i,j}(y)\,d\mu(y).$$

We  denote $\DD_m^{p,B} = \DD_m^{p} \cap \DD_m^B$ and
$$b_m^p(y) = \sum_{i:\,Q_{i,m}\in\DD_m^{p,B} } v_{i,m}(y).$$
Recall also that $\vphi_{y,m}^p(x) = \vphi_{y_i,m}(x)$
if $y\in Q_{i,m}$ with $Q_{i,m}\in \DD_m^{p}$, and $\vphi_{y,m}^p(x)=0$ if there
does not exist any cube of the subfamily $\DD^p_m$ containing $y$.
Then we have
$$U^B_m(x) = \sum_{p=1}^N \int \vphi_{y,m}^p(x)\, b_m^p(y)\,d\mu(y).$$

Now we set $h_m^p = g_m^p + b_m^p$, and we get
$$f(x) = h_0(x) + \sum_{p=1}^N
\sum_{m=1}^\infty \int \vphi_{y,m}^p(x)\, h_m^p(y)\,d\mu(y),$$
with $C\,\vphi_{y,m}\sim y$ for some constant $C>0$, and
$$|h_0| + \sum_{p=1}^N \sum_{m=1}^\infty |h_m^p| \leq C\,A\,\|f\|_*,$$
and the Main Lemma follows, by (g) and \rf{cor2'}.


\subsection{The construction of $g_m$ and $b_m$}

In this subsection we will construct inductively functions $g_m$ and $b_m$
satisfying the properties (a)--(e).
We will check in Subsection
\ref{subsecb} that these functions fulfil (f)--(h) too.

Assume that $g_1,\ldots,g_{m-1}$ and $b_1,\ldots,b_{m-1}$ have been constructed
and they satisfy (a)--(e). 
Let $\Omega_m$ be the set of points $x\in \supp(\mu)$ with $\delta(x,2R_0) > m\,A$ such that
that there exists some $Q \in \DD_m$, $\ell(Q)>0$,
with $Q\ni x$ and $|m_Q f_m|\geq
\frac{3}{4}A$. For each $x\in \Omega_m$,
we consider a doubling cube $S_{x,m}$ centered at $x$
such that $\delta(S_{x,m},2R_0) = mA -\alpha_1-\alpha_2-\alpha_3 \pm \ve_1,$
where $\alpha_3$ is some big constant with $10\alpha_2<\alpha_3\ll A$, whose
precise value will be  fixed below. One has to think that
$S_{x,m}$ is
much bigger than $Q_{x,m}^3$ but much smaller than $Q_{x,m-1}$ (observe that
all these cubes have positive side length).

Now we take a
Besicovich covering of $\Omega_m$ with cubes of type $S_{x,m}$, $x\in\Omega_m$:
$$\Omega_m \subset \bigcup_j S_{j,m},$$
where $S_{j,m}$ stands for $S_{x_{j},m}$, with $x_{j}\in \Omega_m$.
We say that a cube $Q\in \DD_m$ is good (i.e. $Q\in \DD_m^G$) if
$$Q\subset \bigcup_j \frac{3}{2}S_{j,m},$$
and we say that it is bad (i.e. $Q\in \DD_m^B$) if it is not good and
$$Q\subset \bigcup_j  2S_{j,m}.$$
Both good and bad cubes are contained in $\bigcup_j 2S_{j,m}$.
Roughly speaking, the difference between good and bad cubes is that bad 
cubes may be
supported near the boundary of $\bigcup_j 2S_{j,m}$, while the good ones
are far from the boundary.

Now we define $g_m$ and $b_m$:
$$g_m= \sum_{i:\,Q_{i,m}\in \DD_m^G} w_{i,m}\, m_{Q_{i,m}}(f_m),$$
$$b_m= \sum_{i:\,Q_{i,m}\in \DD_m^B} w_{i,m}\, m_{Q_{i,m}}(f_m).$$
Because there is some overlapping
among the cubes in
$\DD_m$, we have used the weights $w_{i,m}$ in the definition of these
functions. However one should think that $g_m$ and $b_m$ are approximations of
the mean of $f$ over the cubes of $\DD_m^G$ and $\DD_m^B$, respectively.

The following remark will be useful.


\begin{claim} \label{cl1.5}
Let $Q_{h,m}\in \DD_m$ be such that either
$g_m\not\equiv0$, $b_m\not\equiv0$ or
$U_m\not\equiv0$ on $Q_{h,m}$. Then there exists some $j$ such that
$\wh{Q}^3_{h,m}\subset 4S_{j,m}$ and so $Q_{h,m}\subset 4S_{j,m}$.
\end{claim}

\begin{proof}
In the first two cases $Q_{h,m}\cap 2S_{j,m}\neq\varnothing$ for some $j$.
In the latter case, by (a) of Lemma \ref{phi} and our construction,
there exists
some $j$ such that $\wh{Q}^3_{h,m} \cap 2S_{j,m}\neq\varnothing$.

So in any case $\wh{Q}^3_{h,m} \cap 2S_{j,m}\neq\varnothing$ for some $j$.
Arguing as in Lemma \ref{mides}, for $\alpha_3$ big enough,
it is easily checked that $\ell(\wh{Q}_{h,m}^3) \leq \ell(S_{j,m})/4$,
and so $\wh{Q}^3_{h,m}\subset 4S_{j,m}$.
\end{proof}


Let us see now that (e) is satisfied.

\begin{claim} \label{claime}
If $Q\in \DD_m$ and $\delta(Q,2R_0)\leq (m-\frac{1}{10})\,A$
(so $\ell(Q)=0$), then
$U_m\equiv g_m\equiv b_m\equiv 0$ on $Q$ and $Q\not\in \DD_m^G\cup \DD_m^B$.
\end{claim}

\begin{proof} Assume that $Q\equiv \{x\}$ and that either
$g_m\not\equiv0$, $b_m\not\equiv0$ or
$U_m\not\equiv0$ on $Q$, or  $Q\in \DD_m^G\cup \DD_m^B$.
By the preceding claim, $Q\subset4S_{j,m}$ for some $j$. Then,
\begin{eqnarray*}
\delta(x,2R_0) & = & \delta(x,4S_{j,m}) + \delta(4S_{j,m},2R_0) \pm\ve_0\\
& \geq & \delta(4S_{j,m},2R_0) - \ve_0 \\
& \geq & \delta(S_{j,m},2R_0) - 8^n\,C_0 - \ve_0 >
\left(m-\frac{1}{10}\right)\,A.
\end{eqnarray*}
\end{proof}

The following estimate will be necessary in many steps of our construction.


\begin{claim} \label{cl1}
Let $Q$ be some cube of the $m$-th generation and $x,y\in 2Q$. Then,
if $g_1,\ldots,g_m$ and $b_1,\ldots,b_m$ satisfy (a), then
$$\sum_{k=1}^m |U_k(x) - U_k(y)|\leq \frac{A}{100}\,\|f\|_*.$$
\end{claim}

We postpone the proof of Claim \ref{cl1} until Subsection \ref{remclaims}.
Let us see that (a) holds.


\begin{claim}  \label{cl4.5}
If $Q\in \DD_m^G\cup \DD_m^B$, then $|m_Q f_m| \leq C_9\,A\,\|f\|_*$. Also,
$|g_m|, |b_m|\leq C_8\, A\,\|f\|_*$.
\end{claim}

\begin{proof}
First we will prove the first statement.
By Claim \ref{claime}, we know that $\delta(Q,2R_0) >(m-\frac{1}{10})\,A$.
Let $R\in \DD_{m-1}$
be such that $Q\cap R\neq \varnothing$. We must have $\ell(R)>0$. Otherwise,
$Q\equiv R$ and
$\delta(R,2R_0) > (m-\frac{1}{10})\,A > (m-1)\,A+\ve_1$, which is not possible.

Since $\ell(Q)\leq \ell(R)/10$, we have
$Q\subset 2R$.
We know $|m_R f_m| \leq A\,\|f\|_*$ because (b) holds for $m-1$.
By Claim \ref{cl1} (for $m-1$ and $R$) we get
\begin{eqnarray*}
|m_Q f_m| & \leq & |m_R f_m| + |m_Q f_m - m_R f_m| \\
& \leq &|m_R f_m| + |m_Q f - m_R f| + \Bigl| m_Q \Bigl( \sum_{k=1}^{m-1} U_k
\Bigr)  - m_R \Bigl( \sum_{k=1}^{m-1} U_k\Bigr) \Bigr|  \\
& \leq & C\,A\,\|f\|_* + |m_Q f - m_R f|.
\end{eqnarray*}
The term $|m_Q f - m_R f|$ is also bounded above by $C\,A\,\|f\|_*$
because $Q$ and $R$ are doubling, $f\in\rbmo(\mu)$, and
it is easily checked that $\delta(Q,R) \leq C\, A$.

The estimates on $g_m$ and $b_m$ follow from from the definition of these
functions and the estimate
$|m_Q f_m| \leq C_9\,A\,\|f\|_*$ for $Q\in\DD_m^G\cup\DD_m^B$.
\end{proof}

Let us prove (d) now.


\begin{claim}
If $Q\in \DD_m$ and $|m_Q f_m|\leq \dfrac{8}{20}\,A\,\|f\|_*$, then
$U_m\equiv 0$ and $g_m\equiv b_m\equiv 0$ on $Q$.
\end{claim}

\begin{proof}
Suppose that $Q\equiv Q_{h,m}\in \DD_m$ is such that either
$g_m\not\equiv0$, $b_m\not\equiv0$ or
$U_m\not\equiv0$ on $Q_{h,m}$.
By Claim \ref{cl1.5} we have $Q_{h,m}\subset 4S_{j,m}$ for some $j$. By
construction, the center of $S_{j,m}$ belongs to some cube $Q_{i,m}$ with
$|m_{Q_{i,m}} f_m|\geq \frac{3}{4} \,A\,\|f\|_*$.
It is easily seen that $\delta(Q_{h,m},4S_{j,m}),\,\delta(Q_{i,m},4S_{j,m})
\leq C' + \alpha_1 + \alpha_2 + \alpha_3$. Thus
$$|m_{Q_{i,m}} f- m_{Q_{h,m}} f| \leq (C'' + 2\alpha_1 + 2\alpha_2 + 2\alpha_3)
\,\|f\|_*.$$
Since $Q_{i,m}$ and $Q_{h,m}$ are contained in a common cube of
the generation $m-1$, by Claim \ref{cl1} we get
\begin{eqnarray*}
|m_{Q_{i,m}} f_m- m_{Q_{h,m}} f_m| & \leq & |m_{Q_{i,m}} f- m_{Q_{h,m}} f|
\\ && \mbox{}+
\Bigl| m_{Q_{i,m}} \Bigl( \sum_{k=1}^{m-1} U_k
\Bigr)  - m_{Q_{h,m}} \Bigl( \sum_{k=1}^{m-1} U_k\Bigr) \Bigr| \\
& \leq & (C'' + 2\alpha_1 + 2 \alpha_2 + 2\alpha_3 + A/100)\, \|f\|_* \\
& \leq & \frac{1}{10}\,A\,\|f\|_*,
\end{eqnarray*}
and so
$$|m_{Q_{h,m}} f_m| \geq \left(\frac{3}{4} - \frac{1}{10}\right)\,A\,\|f\|_* >
\frac{8}{20}\,A\,\|f\|_*.$$
\end{proof}

The statement (c) is a consequence of the fact that if
$Q\in \DD_m^G$, then $Q$ is far  from the boundary of $\bigcup_j2 S_{j,m}$.
Then $U_m$ is very close to $m_Q f_m$ on $Q$,
since we only integrate over cubes of
$\DD^G_m\cup\DD_m^B$ in order to obtain $U_m(x)$ for $x\in Q$.
On the other hand, if $Q\in\DD^B_m$, this argument does not work because
$Q$ may be near the boundary of $\bigcup_j2 S_{j,m}$, and so it may
happen that we integrate on some cubes from
$\DD_m\setminus (\DD_m^G\cup\DD_m^B)$ for obtaining $U_m(x)$, $x\in Q$.

Let us see (c) in detail.


\begin{claim} \label{cl2}
If $Q\in \DD_m^G$ and $\ell(Q)>0$, then $|m_Q f_{m+1}|
\leq \frac{7}{20}A \|f\|_*$.
\end{claim}

\begin{proof} Consider $Q_{i,m} \in \DD_m^G$. We want to see that
$U_m$ is very close to $m_{Q_{i,m}} f_m$ on this cube. By (a) of Lemma \ref{phi}
we have to deal with the cube $\wh{Q}_{i,m}^3$.

Let us see that if $P\in\DD_m$ is such that
$P\cap\wh{Q}_{i,m}^3\neq\varnothing$, then $P\in \DD_m^G\cup \DD_m^B$.
Notice that $P\subset {\QH}_{i,m}$. Now, by
the definition of good cubes, there exists
some $j$ such that $Q_{i,m}\cap \frac{3}{2}S_{j,m}\neq\varnothing$, which
implies ${\QH}_{i,m}\cap \frac{3}{2}S_{j,m}\neq\varnothing$.
For $\alpha_3$ big enough, we have $\ell({\QH}_{i,m})\ll \ell(S_{j,m})$, and
then ${\QH}_{i,m}\subset 2S_{j,m}$. So $P\in \DD_m^G\cup \DD_m^B$.

Let us estimate the term
$$\sup_{y\in \wh{Q}_{i,m}^3} |(g_m(y)+b_m(y)) - m_{Q_{i,m}} f_m|.$$
Recall that
$$g_m(y) + b_m(y)= \sum_{h:\,Q_{h,m}\in \DD_m^G\cup\DD_m^B}
 w_{h,m}(y)\, m_{Q_{h,m}}f_m.$$
By the arguments above, if $y\in \wh{Q}_{i,m}^3$ and $w_{h,m}(y)\neq0$, then
$Q_{h,m}$ has been chosen for supporting $g_m$ or $b_m$, i.e.
$Q_{h,m}\in \DD_m^G\cup\DD_m^B$. Then,
$$g_m(y) + b_m(y) -  m_{Q_{i,m}} f_m =
\sum_{h:\,Q_{h,m}\in \DD_m}
w_{h,m}(y)\, (m_{Q_{h,m}}f_m - m_{Q_{i,m}} f_m).$$
By Claim \ref{cl1} we obtain
\begin{eqnarray*}
|m_{Q_{h,m}}f_m - m_{Q_{i,m}} f_m| & \leq & \frac{1}{100}\,A \|f\|_* +
|m_{Q_{h,m}}f - m_{Q_{i,m}} f| \\
& \leq &
\left(\frac{1}{100}\,A + C + 2\,\delta(Q_{h,m},Q_{i,m}) \right)\,\|f\|_* \\
& \leq & \frac{1}{50} \,A\,\|f\|_*
\end{eqnarray*}
(we have used that $\delta(Q_{h,m},Q_{i,m}) \leq C$, with
$C$ depending on $\alpha_1,\,\alpha_2$). Then we get
\begin{equation}  \label{suppp}
|g_m(y) + b_m(y) -  m_{Q_{i,m}} f_m |\leq \frac{1}{50} \,A\,\|f\|_*.
\end{equation}

For $x\in Q_{i,m}$, we have
\begin{eqnarray} \label{sup2}
|U_m(x) - m_{Q_{i,m}} f_m| & \leq &
\left|U_m(x) - m_{Q_{i,m}} f_m\, \int \vphi_{y,m}(x)\,d\mu(y) \right| \nonumber
\\ && \mbox{}+ |m_{Q_{i,m}} f_m|\,
\left|1- \int \vphi_{y,m}(x)\,d\mu(y)\right|.
\end{eqnarray}
Let us estimate the first term on the right hand side. By \rf{suppp} and
\rf{convo1} we obtain
\begin{eqnarray*}
\lefteqn{\left|U_m(x) - m_{Q_{i,m}} f_m\, \int \vphi_{y,m}(x)\,d\mu(y) \right|
} \\ & = &
\left|\int_{\wh{Q}_{i,m}^3}
 \vphi_{y,m}(x) \,(g_m(y)+b_m(y)- m_{Q_{i,m}} f_m)\,d\mu(y) \right|\\
& \leq & (1+\ve_3)\, \frac{1}{50} \,A\,\|f\|_* .
\end{eqnarray*}
On the other hand, by \rf{convo1}, \rf{convo2} and Claim \ref{cl4.5},
the second term on the right hand side of \rf{sup2} is bounded
above by $\ve_3\,C_8\,A\,\|f\|_*$. Thus we have
$$|m_{Q_{i,m}} f_{m+1}| \leq \left((1+\ve_3)\, \frac{1}{50} +
\,\ve_3\,C_8\right)
 \,A\,\|f\|_* \leq  \frac{7}{20}\,A\,\|f\|_*,$$
if we choose $\ve_3$ small enough.
\end{proof}

Now we are going to show that (b) also holds.


\begin{claim} \label{cl3}
If $Q\in \DD_m$ and $\ell(Q)>0$, then $|m_Q f_{m+1}| \leq A \|f\|_*$.
\end{claim}

\begin{proof}
If $Q\in \DD_m^G$, we have already seen that $|m_Q f_{m+1}|
\leq
\frac{7}{20}\,A\,\|f\|_*$.

If $Q\in \DD_m\setminus\DD_m^G$, then $Q\cap \bigcup_j S_{j,m}=\varnothing$
(because $\ell(Q) \ll \ell(S_{j,m})$ and $Q\not\subset
\bigcup_j \frac{3}{2}S_{j,m}$). By construction, we have
\begin{equation}  \label{defyy0}
|m_Q f_m|\leq \frac{3}{4}\,A\,\|f\|_*.
\end{equation}
If $U_m\equiv0$ on $Q$, then
$|m_Q f_{m+1}|=|m_Q f_m|\leq \frac{3}{4}\,A\,\|f\|_*$.

Now we consider the case
$Q\equiv Q_{h,m}\cap \bigcup_j S_{j,m}=\varnothing$ such that
$U_m\not\equiv0$ on $Q$.
By Claim \ref{cl1.5} there exists some $j$ with
$\wh{Q}_{h,m}^3\subset 4S_{j,m}$.
Recall that by (a) of Lemma \ref{phi}, if $x\in Q_{h,m}$, we have
$$U_m(x) = \int_{\wh{Q}_{h,m}^3} \vphi_{y,m}(x)\, (g_m(y) + b_m(y))\, d\mu(y).$$
So if $\vphi_{y,m}(x)\neq 0$ and $y\in Q_{i,m}$, we have $Q_{i,m} \cap
\wh{Q}_{h,m}^3\neq \varnothing$. Therefore, $Q_{i,m} \subset
{\QH}_{h,m}$. Then,
$$\delta(Q_{i,m},Q_{h,m}) \leq C + \delta(Q_{i,m},{\QH}_{h,m})
+ \delta(Q_{h,m},{\QH}_{h,m}) \leq  C + 2\alpha_1 + 2\alpha_2 \leq
\frac{A}{400}.$$
Therefore, $|m_{Q_{i,m}} f - m_{Q_{h,m}}f| \leq \frac{A}{100}\,\|f\|_*$.
By Claim 1 we get
\begin{eqnarray}  \label{defyy}
|m_{Q_{i,m}} f_m- m_{Q_{h,m}} f_m| & \leq & |m_{Q_{i,m}} f- m_{Q_{h,m}} f|
\nonumber \\ && \mbox{}+
\Bigl| m_{Q_{i,m}} \Bigl( \sum_{k=1}^{m-1} U_k
\Bigr)  - m_{Q_{h,m}} \Bigl( \sum_{k=1}^{m-1} U_k\Bigr) \Bigr| \nonumber \\
& \leq & \frac{1}{10}\,A\,\|f\|_*.
\end{eqnarray}
Recall also that, by (d),
\begin{equation}  \label{defyy2}
|m_{Q_{h,m}} f_m| \geq \frac{8}{20}\,A\,\|f\|_*.
\end{equation}
From the definition of $g_m, b_m$ and \rf{defyy}, \rf{defyy2}, we derive that
$m_{Q_{h,m}} f_m$ and $U_m(x)$ have the same sign.

On the other hand, from \rf{defyy0} and \rf{defyy} we get
$$|m_{Q_{i,m}} f_m|\leq \frac{34}{40}\,A\,\|f\|_*.$$
So by the definition of $g_m$ anb $b_m$ we have
$$\|g_m + b_m\|_{L^\infty(\mu)} \leq
\frac{34}{40}\,A\,\|f\|_*,$$
and by \rf{convo1} we obtain
\begin{equation}  \label{deffi}
|U_m(x)| \leq  \frac{34}{40}\,A\,\|f\|_*\,\int \vphi_{y,m}(x)\,d\mu(y)
\leq  (1+\ve_3)\,\frac{34}{40}\,A\,\|f\|_* \leq A\,\|f\|_*
\end{equation}
(assuming $\ve_3$ small enough).
By \rf{defyy0}, \rf{deffi} and since $m_{Q_{h,m}} f_m$ and $U_m(x)$ have the
same sign, (b) holds also in this case.
\end{proof}

Therefore, (a)--(e) are satisfied.


\subsection{Proof of (f), (g) and (h)} \label{subsecb}

The statement (f) is a direct consequence of the following.

\begin{claim}
If $\delta(x,2R_0)<\infty$,
and if $Q = \{x\}\in\DD_m$  (i.e. $\ell(Q)=0$),
then $h_0(x) = f_{m+1}(x)$ and $|h_0(x)| \leq C_{9}\,A\,\|f\|_*.$
\end{claim}

\begin{proof}
Take $m$ such that $(m-1)\,A< \delta(x,2R_0) \leq m\,A$. By (e) we get
$U_{m+k}(x)=0$ for $k\geq1$. Therefore,
$f_{m+1}(x) = f_{m+2}(x) = \cdots = h_0(x).$
By (a) we have
$$|f_{m+1}(x)|\leq |f_m(x)| + |U_m(x)| \leq |f_m(x)| +
2\,C_8\,(1+\ve_3)\,A\,\|f\|_*.$$
So we only have to estimate $|f_m(x)|$.

Take $Q_{i,m-1}\in \DD_{m-1}$ with $x\in Q_{i,m-1}$.
Since $\ell(Q_{i,m-1})>0$, by (b) we have $|m_{Q_{i,m-1}}f_m| \leq A\,\|f\|_*$.
Applying Claim \ref{cl1} we get
\begin{eqnarray*}
|m_{Q_{i,m-1}}f_m - f_m(x)| & \leq & |m_{Q_{i,m-1}}f - f(x)| +
\frac{A}{100}\,\|f\|_* \\
& \leq & C\,\left(1+\delta(x,Q_{i,m-1}) +\frac{A}{100}\right)\,\|f\|_*.
\end{eqnarray*}
It is easily checked that $\delta(x,Q_{i,m-1}) \leq A + \ve_0 + \ve_1$.
Then we get $|f_m(x)|\leq C\,A\,\|f\|_*$.
\end{proof}


Now we turn our attention to (g).
Given some good cube $Q_{i,m}\in \DD_m^G$ with $\ell(Q_{i,m})>0$, we denote
$$Z_{i,m}:= Z(Q_{i,m},A\,\|f\|_*/30)$$
(see Definition \ref{zq}; roughly speaking $Z_{i,m}$ is the
part of $Q_{i,m}$
where $f$ does not oscillate too much with respect to $m_{Q_{i,m}} f$).
If $Q_{i,m}\in \DD_m^G$ and $\ell(Q_{i,m})=0$, we set $Z_{i,m}=Q_{i,m}$.
The set $Z_{i,m}$ has a very nice property:


\begin{claim} \label{cl4}
Let $k>m$ and $Q_{i,m}\in \DD_m^G$. If $P\in\DD_k$ is such that $P\cap
Z_{i,m}\neq\varnothing$, then $g_k\equiv b_k\equiv 0$ on $P$ and $P\not\in
\DD_k^G\cup\DD_k^B$.
\end{claim}

\begin{proof} Consider first the case $\ell(Q_{i,m}) = 0$.
If $P\in\DD_k$ is such that $P\cap
Q_{i,m}\neq\varnothing$, then $\ell(P)\leq \ell(Q_{i,m})/10=0$ and so
$P\equiv Q_{i,m}$. Therefore,
$$\delta(P,2R_0) \leq m\,A \leq \left(k-\frac{1}{10}\right)\,A.$$
By (e), we get $b_k\equiv g_k\equiv 0$ on $P$.

\vv
Assume now $\ell(Q_{i,m})>0$.
Let $x\in P\cap Z_{i,m}$. From the definition of $Z_{i,m}$, we have
\begin{equation} \label{cvb1}
|m_{Q_{i,m}} f - m_S f| \leq \frac{A}{30}\,\|f\|_*
\end{equation}
for any $S\in\DD_{m+j}$, $j\geq1$, with $x\in S$. Also, by Claim \ref{cl2}
we have
$$|m_{Q_{i,m}} f_{m+1}| \leq \frac{7}{20}A\,\|f\|_*.$$
Consider now 
$P_{m+1}\in \DD_{m+1}$ with $x\in P_{m+1}$. Observe that $\ell(P_{m+1})\leq
\ell(Q_{i,m})/10$ and $P_{m+1}\subset 2Q_{i,m}$.
We have
\begin{eqnarray*}
|m_{P_{m+1}} f_{m+1}| & \leq & |m_{Q_{i,m}} f_{m+1}| + |m_{Q_{i,m}} f_{m+1}
- m_{P_{m+1}} f_{m+1}| \\
& \leq & \frac{7}{20}A\,\|f\|_* +  |m_{Q_{i,m}} f - m_{P_{m+1}} f| \\ &&\mbox{}+
\Bigl| m_{Q_{i,m}} \Bigl( \sum_{k=1}^{m} U_k
\Bigr)  - m_{P_{m+1}} \Bigl( \sum_{k=1}^{m} U_k\Bigr) \Bigr|.
\end{eqnarray*}
By \rf{cvb1} and Claim \ref{cl1} we obtain $|m_{P_{m+1}}
f_{m+1}|\leq\frac{8}{20}A\,\|f\|_*$. By (d), on $P_{m+1}$ we have
$g_{m+1} \equiv b_{m+1} \equiv0$ and also $U_{m+1}\equiv 0$.
Thus,
$$f_{m+2} \equiv f_{m+1}$$
on any cube $P_{m+1}\in\DD_{m+1}$ containing $x$.

Take now $P_{m+2}\in \DD_{m+2}$ with $x\in P_{m+2}$. On this cube
$f_{m+2} \equiv f_{m+1}$, and then we have
\begin{eqnarray*}
|m_{P_{m+2}} f_{m+2}| & \leq & |m_{Q_{i,m}} f_{m+1}| + |m_{Q_{i,m}} f_{m+1}
- m_{P_{m+2}} f_{m+1}| \\
& \leq & \frac{7}{20}A\,\|f\|_* +  |m_{Q_{i,m}} f - m_{P_{m+2}} f| \\
&& \mbox{} + \Bigl| m_{Q_{i,m}} \Bigl( \sum_{k=1}^{m} U_k
\Bigr)  - m_{P_{m+2}} \Bigl( \sum_{k=1}^{m} U_k\Bigr) \Bigr|.
\end{eqnarray*}
Again by (d), we get $g_{m+2} \equiv b_{m+2} \equiv U_{m+2}
\equiv0$ on $P_{m+2}$. Thus, $f_{m+3} = f_{m+1}$ on $P_{m+2}$.

Going on, we will obtain $g_{m+j} \equiv b_{m+j} \equiv U_{m+j}
\equiv0$ for all $j\geq1$ on any cube $P_{m+j}\in \DD_{m+j}$ containing $x$.
\end{proof}

As a consequence of Claim \ref{cl4}, $Z_{i,m}$ is a good place for
supporting $g_m$. If, for each
$m$, $g_m$ were supported on $\bigcup_i Z_{i,m}$, then the supports of
$g_m$, $m\geq1$, would be disjoint for different $m$'s. This is the idea that
Carleson used in \cite{Carleson}.

So we are going to make some ``corrections'' according to this argument.
We have
$$U_m^G(x) = \sum_{i\in I_m} \vphi_{y_i,m}(x) \int w_{i,m}(y)\, g_m(y)\, d\mu(y).$$
For each $Q_{i,m}$ with $\ell(Q_{i,m})>0$ we set
$$u_{i,m}(y) = \int w_{i,m}\, g_m\, d\mu\,\cdot
\frac{\chi_{Z_{i,m}}(y)}{\mu(Z_{i,m})}.$$
If $\ell(Q_{i,m})=0$, we set $u_{i,m}(y) = w_{i,m}(y)\, g_m(y) \equiv g_m(y)$
(we do not change anything in this case).
Then $U_m^G$ can be written as
$$U_m^G(x) = \sum_{i\in I_m} \vphi_{y_i,m}(x) \int u_{i,m}(y)\, d\mu(y).$$
As in the case of $U_m^B$ in Subsection \ref{correccio}, if we set
$\DD^G_m = \DD_m^{1,G} \cup \cdots \cup\DD_m^{N,G}$ where each subfamily
$\DD_m^{p,G}$ is disjoint, we can write $U_m^G$ in the following way:
$$U^G_m(x) = \sum_{p=1}^N \int \vphi_{y,m}^p(x)\, g_m^p(y)\,d\mu(y)$$
with
$$g_m^p(y) = \sum_{i:\,Q_{i,m}\in\DD_m^{p,G} } u_{i,m}(y)$$
and
$$\vphi_{y,m}^p(x) = \vphi_{y_i,m}(x)$$
if $y\in Q_{i,m}$ and $Q_{i,m}\in \DD_m^{p}$. 

By Proposition \ref{JN2}, if $A$ is big enough we have $\mu(Z_{i,m})\geq
\mu(Q_{i,m})/2$ (if $\ell(Q_{i,m})>0$). Then
it easily checked that $\|u_{i,m}\|_{L^\infty(\mu)}\leq 2\,
\|g_m\|_{L^\infty(\mu)}$
for all $i$. Thus, from (a), (g.2) follows.
Moreover, because of Claim \ref{cl4}, (g.3) also holds.

\vspace{2mm}
One of the differences between our construction and Carleson's one
is that, because of the regularity
of Lebesgue measure, Carleson can treat the bad cubes in a way very similar to
the way for the good ones.
We have not been able to operate as Carleson. However, as it has
been shown in Subsection \ref{correccio}, the packing condition \rf{pack}
is also a good solution. Let us prove that this condition is satisfied.


\begin{claim} \label{cl5}
For any $R\in\DD_m$ with $\ell(R)>0$, the bad cubes satisfy the packing condition
$$\sum_{\substack{ Q:\,Q\cap R\neq\varnothing \\Q\in D^B_k,\,k> m}}
\mu(Q) \leq C\,\mu(R).$$
\end{claim}

\begin{proof}
Let $k>m$ be fixed. We are going to estimate the sum
$$\sum_{\substack{ Q:\,Q\cap R\neq\varnothing \\Q\in D^B_k}}
\mu(Q).$$
Let $Q\in \DD^B_k$ be such that $Q\cap R\neq\varnothing$.
Since $Q$ is a bad cube, there exists some
$j$ such that $2S_{j,k} \cap Q\neq \varnothing$.
Then we have $Q\subset 4S_{j,k}$.
Since $A\gg \alpha_1+\alpha_2+\alpha_3$ and $4S_{j,k}\cap R\neq\varnothing$, we
get $\ell(S_{j,k})\leq \ell(R)/20$, and so $4S_{j,k}\subset 2R$.

By the finite overlapping of the cubes $Q$ in $\DD_k$, we have
\begin{multline*}
\sum_{\substack{ Q:\,Q\cap R\neq\varnothing \\Q\in D^B_k}} \mu(Q)
\,\leq \, C\, \mu\biggl( \bigcup_{j:\, S_{j,k}\subset 2R} 2S_{j,k} \biggr)  \\
  \leq\,  C\, \sum_{j:\, S_{j,k}\subset 2R} \mu(2S_{j,k})
\, \leq \, C\, \sum_{j:\, S_{j,k}\subset 2R} \mu(S_{j,k}).
\end{multline*}
Now, from the construction of $g_k^p$, it is easy to check that $\mu(S_{j,k})
\leq C\,\mu\Bigl(S_{j,k} \cap
\bigl\{\sum_{p=1}^N |g_k^p|\neq 0\bigr\}\Bigr)$.
This fact and the bounded overlapping of the cubes $S_{j,k}$ give
$$\sum_{\substack{ Q:\,Q\cap R\neq\varnothing \\Q\in D^B_k}} \mu(Q)
\leq C\,  \mu\Bigl(2R \cap
\Bigl\{\sum_{p=1}^N |g_k^p|\neq 0\Bigr\}\Bigr).$$
Summing over $k>m$, as the supports of the functions $g_k^p$ are disjoint for
different $k$'s, we obtain
$$\sum_{\substack{ Q:\,Q\cap R\neq\varnothing \\Q\in D^B_k,\,k> m}}
\mu(Q) \leq C\ \sum_{k>m} \mu\Bigl(2R \cap
\Bigl\{\sum_{p=1}^N |g_k^p|\neq 0\Bigr\}\Bigr) \leq C\,\mu(2R) \leq
C\,\mu(R).$$
\end{proof}


\subsection{Proof of Claim \ref{cl1}}  \label{remclaims}

We only need to check that
$$\sum_{k=1}^m
C_8\, A\,\int |\vphi_{z,k}(x) - \vphi_{z,k}(y)|\, d\mu(z) \leq\frac{A}{100}.$$
Let $x_0\in \supp(\mu)$ be such that $x,y\in 2Q_{x_0,m}$.
Obviously, we can assume $\ell(Q_{x_0,m})>0$.
For each $k\leq m$ we set
$$\int |\vphi_{z,k}(x) - \vphi_{z,k}(y)|\, d\mu(z) =
\int_{\R^d\setminus \check{Q}_{x_0,k}^1} +
\int_{\check{Q}_{x_0,k}^1}= I_{1,k} + I_{2,k}.$$

Let us estimate the integrals $I_{1,k}$. Notice that
if $x,y\in 2Q_{x_0,m}$, then $x,y\in 2Q_{x_0,k} \subset \frac{1}{2}
\check{Q}_{x_0,k}^1$. Thus $|x-z|\approx|y-z|\approx|x_0-z|$ for $z\in
\R^d\setminus \check{Q}_{x_0,k}^1$. So by (d) of Lemma \ref{phi} we have
\begin{eqnarray}  \label{ii1}
I_{1,k} & \leq & C\,\alpha_2^{-1}\,\int_{\R^d\setminus \check{Q}_{x_0,k}^1}
\frac{|x-y|}{|x-z|^{n+1}}\,d\mu(z) \nonumber \\
& \leq & C\,\alpha_2^{-1}\,\frac{\ell(Q_{x_0,m})}{\ell(\check{Q}_{x_0,k}^1)}.
\end{eqnarray}
In case $k>m$, by Lemma \ref{mides2} we get
$$I_{1,k} \leq C \,\alpha_2^{-1}\,\frac{\ell(Q_{x_0,m})}{\ell(Q_{x_0,k})} \leq
C_{13}\,\alpha_2^{-1}\,
2^{-\gamma\,(m-k)\,A}.$$
Therefore,
\begin{equation} \label{uu1}
C_8\, A\,\sum_{k=1}^m I_{1,k} \leq C_8\, \alpha_2^{-1}\,A\,
\sum_{k=1}^{m-1}2^{-\gamma\,(m-k)\,A} + C_8\,C_{13}\,\alpha_2^{-1}\,A\,
\frac{\ell(Q_{x_0,m})}{\ell(\check{Q}_{x_0,m}^1)}.
\end{equation}
The first sum on the right hand side is $\leq
C\,\alpha_2^{-1}\,A\,2^{-\gamma\,A}$, and
for $A$ big enough and $\alpha_2>1$ is $\leq 1 \leq A/400$.
The second term on the right hand side is also $\leq A/400$ if we choose
$\alpha_2$ big enough (or $\alpha_1$ big enough since then
$\ell(\check{Q}_{x_0,m})\gg \ell(\check{Q}_{x_0,m}^1)$).
Thus
$$C_8\, A\,\sum_{k=1}^m I_{1,k} \leq \frac{A}{200}.$$


We consider now the integrals $I_{2,k}$. By Lemma \ref{phi},
$$|\vphi'(u)|\leq
C\,\frac{ \alpha_2^{-1}}{\ell(\check{Q}_{x_0,k}^1)^{n+1}}$$
for all $u\in Q_{x_0,k}$. Therefore,
$$I_{2,k} \leq C\,\alpha_2^{-1}\,\int_{\check{Q}_{x_0,k}^1}
\frac{|x-y|}{\ell(\check{Q}_{x_0,k}^1)^{n+1}} \,d\mu(z) \leq C\,\alpha_2^{-1}\,
\frac{\ell(Q_{x_0,m})}{\ell(\check{Q}_{x_0,k}^1)}.$$
This is the same estimate that we have obtained for $I_{1,k}$ in \rf{ii1}, and
then we also have
$$C_8\, A\,\sum_{k=1}^m I_{2,k} \leq \frac{A}{200},$$
if we choose $A$ and $\alpha_2$ (or $\alpha_1$) big enough.
\qed


\section{Appendix}  \label{secappen}

In this section we will prove the following result, which is used in Section
\ref{secdual}
to show that Theorem \ref{maxi} follows from the Main Lemma.

\begin{lemma}  \label{final}
Consider $f\in L^1(\mu)$ with $\int f\,d\mu=0$ and $M_\Phi f \in L^1(\mu)$.
Then there exists
a sequence of functions $f_k$, $k\geq1$, bounded with compact support such that
$\int f_k\,d\mu=0$, $f_k\to f$ in $L^1(\mu)$ and $\|M_\Phi
(f-f_k)\|_{L^1(\mu)}\to0$.
\end{lemma}

So if we consider the space
$$H^1_\Phi(\mu) = \Bigl\{ f\in L^1(\mu):\, {\textstyle \int
f\,d\mu=0},\,M_\Phi f\in L^1(\mu)\Bigr\},$$
with norm $\|f\|_{H^1_\Phi(\mu)} = \|f\|_{L^1(\mu)} + \|M_\Phi
f\|_{L^1(\mu)}$, then Lemma \ref{final} asserts that functions in
$H^1_\Phi(\mu)$ which are bounded and have compact support are dense in
$H^1_\Phi(\mu)$. In particular, $H^1_\Phi(\mu)\cap \hbm$ is
dense in $H^1_\Phi(\mu)$.

In this section we will assume that the center of any cube $Q$ may be any
point of $\R^d$, not necessarily belonging to $\supp(\mu)$.
As in the previous sections, the sides of the cubes are
parallel to the axes and they are closed.

Let us introduce some additional notation. For $\rho>1$, we set
$$M_{(\rho)} f(x) = \sup_{Q\ni x} \frac{1}{\mu(\rho Q)}\int_Q |f|\,d\mu.$$
This non centered maximal operator is bounded above by the operator
defined as
$$M^{(\rho)} f(x) = \sup_{\rho^{-1}Q\ni x} \frac{1}{\mu(Q)}\int_Q
|f|\,d\mu.$$
This is the version of the Hardy-Littlewood operator that one
obtains taking supremums over cubes $Q$ which may be non centered at $x$
but such that $x\in \rho^{-1} Q$. Recall that since $0<\rho^{-1}<1$, one can
apply Besicovich's Covering Theorem and then one gets that $M^{(\rho)}$ is of
weak type $(1,1)$ and bounded in $L^p(\mu)$, $p\in(1,\infty]$. As a
consequence, $M_{(\rho)}$ is also of
weak type $(1,1)$ and bounded in $L^p(\mu)$, $p\in(1,\infty]$

\begin{rem}[Whitney covering] \label{witney}
Let $\Omega\subset\R^d$ be open, $\Omega\neq\R^d$. Then $\Omega$ can be
decomposed as $\Omega=\bigcup_{i\in I} Q_i$, where $Q_i,\,i\in I,$ are cubes
with disjoint interiors, with $20Q_i\subset\Omega$ and
such that, for some constants $\beta>20$ and $D\geq1$,
$\beta\,Q_k\cap \Omega^c\neq\varnothing$ and
for each cube $Q_k$ there are at most $D$ cubes $Q_i$ with
$10Q_k\cap10Q_i\neq\varnothing$ (in particular, the family of cubes
$\{10Q_i\}_{i\in I}$ has finite overlapping).
\end{rem}

\vv
In \cite{Tolsa3} a decomposition of Calder\'on-Zygmund type
adapted for non doubling measures was introduced. This
decomposition was used to prove an interpolation theorem between
$(H^1_{atb}(\mu),L^1(\mu))$ and $(L^\infty(\mu),\rbmo(\mu))$. In \cite{Tolsa5} it
was shown
that this decomposition was also useful for proving that CZO's bounded in
$L^2(\mu)$ are of weak type $(1,1)$ too, as in the doubling case (this result
had been proved previously in \cite{NTV2} using different techniques).
To prove Lemma \ref{final} we will use the following variant of the
Calder\'on-Zygmund decomposition of \cite{Tolsa3}.

\begin{lemma}
Let $f\in L^1(\mu)$ with
$\int \!f\,d\mu=0$ and $M_\Phi f\in L^1(\mu)$. For any $\lambda>0$, let
$\Omega_\lambda = \bigl\{x\in\R^d: \,M_{(2)}f(x)>\lambda\bigr\}$.
Then $\Omega_\lambda$ is open and $|f| \leq 2^{d+1}\,\lambda$
$\mu$-a.eq. in $\R^d\setminus\Omega_\lambda$. Moreover, if we consider a
Whitney decomposition of $\Omega_\lambda$ into cubes $Q_i$ (as in Remark
\ref{witney}), then we have:

\begin{itemize}

\item[(a)] For each $i$ there exists a function $w_i\in \CC^\infty(\R^d)$ with
$\supp(w_i)\subset\frac{3}{2}Q_i$, $0\leq w_i\leq 1$, $\|w_i'\|_\infty\leq
C\,\ell(Q_i)^{-1}$ such that
$\sum_i w_i(x) = 1$ if $x\in \Omega_\lambda$.

\item[(b)] For each $i$, let $R_i$ be the smallest $(6,6^{n+1})$-doubling cube
of the form $6^k Q_i$, $k\geq1$, with $R_i\cap \Omega_\lambda^c\ne\varnothing$.
Then there exists a family of functions
$\alpha_i$ with $\supp(\alpha_i)\subset R_i$ satisfying
\begin{equation}  \label{cc4}
\int \alpha_i \,d\mu = \int f\,w_i\,d\mu,
\end{equation}
\begin{equation}  \label{cc4.5}
\|\alpha_i\|_{L^\infty(\mu)} \,\mu(R_i) \leq C\,\|\alpha_i\|_{L^1(\mu)}
\end{equation}
and
\begin{equation}  \label{cc5}
\sum_i |\alpha_i| \leq B\,\lambda
\end{equation}
(where $B$ is some constant).

\item[(c)] $f$ can be written as $f=g+b$, with
$$g = f\,\Bigl(1-\sum_i w_i\Bigr) + \sum_i \alpha_i$$
and
$$b = \sum_i (f\,w_i - \alpha_i),$$
and then $\|g\|_{L^\infty(\mu)}\leq C\,\lambda$ and
$\supp(b)\subset\Omega_\lambda$.

\end{itemize}
\end{lemma}

\begin{proof}
The set $\Omega_\lambda$ is open because $M_{(2)}$ is lower semicontinuous.
Since for $\mu$-a.e. $x\in\R^d$ there exists a sequence of
$(2,2^{d+1})$-doubling cubes centered at $x$ with side length tending to zero,
it follows
that for $\mu$-a.e. $x\in\R^d$ such that $|f(x)|>2^{d+1}\lambda$ there exists
some $(2,2^{d+1})$-doubling cube $Q$ with $\int_Q |f|\,d\mu/\mu(Q)>
2^{d+1}\lambda$ and so $M_{(2)}f(x)>\lambda$.

The existence of the functions $w_i$ of (a) is a standard known fact.
The assertion (c) follows from the other statements in the lemma. So the
only question left is the statement (b).

Notice that, since $R_i\cap \Omega_\lambda^c\neq \varnothing$, we have
\begin{equation}  \label{cc2}
\int_{R_i} |f|\,d\mu \leq \lambda\,\mu(2R_i)
\end{equation}
for each $i$.

To construct the functions $\alpha_j$ we would like to start by the smallest
cube $R_i$, and go on with the bigger cubes $R_j$ following an order of non
decreasing
sizes. Since in general there does not exist a cube $R_i$ with minimal side
length in the family $\{R_i\}_{i=1}^\infty$, 
we will have to modify a little the argument. 
For each fixed $N$ we will construct functions $\alpha_i^N$,
$1\leq i\leq N$, with $\supp(\alpha_i^N)\subset R_i$, satisfying \rf{cc4}, 
\rf{cc4.5}
and \rf{cc5}. Finally, applying weak limits when $N\to\infty$, we will get
the functions $\alpha_i$.

The functions $\alpha_i^N$ that we will construct will be of the form
$\alpha_i^N
=a_i^N\,\chi_{A_i^N}$, with $a_i^N\in\R$ and $A_i^N\subset R_i$.
To avoid a complicate notation, suppose that the cubes $R_i$, $1\leq i \leq N$,
satisify $\ell(R_i)\leq \ell(R_{i+1})$ (we can assume this because we are taking a
finite number of cubes).
We set $A_1^N=R_1$ and
$$\alpha_1^N = a_1^N\,\chi_{R_1},$$
where the constant $a_1^N$ is chosen so that $\int_{Q_1}f\,w_1\,d\mu=
\int \alpha_1^N\,d\mu$.

Suppose that $\alpha_1^N,\, \alpha_{2}^N
\ldots,\alpha_{k-1}^N$ (for some $k\leq N$) have been constructed,
satisfy \rf{cc4} and $\sum_{i=1}^{k-1} |\alpha_i|\leq
B\,\lambda,$  where $B$ is some constant (which will be fixed below).

Let $R_{s_1},\ldots,R_{s_m}$ be the subfamily of cubes
$R_i$, $1\leq i\leq k-1$, such that $R_{s_j}\cap R_k \neq \varnothing$.
As $l(R_{s_j}) \leq l(R_k)$ (because of the non decreasing sizes of $R_i$),
we have $R_{s_j} \subset 3R_k$. Taking into account that for $i=1,\ldots,k-1$
$$\int |\alpha_i^N|\,d\mu \leq \int |f\,w_i|\,d\mu$$
by \rf{cc4}, and using that $R_k$ is $(6,6^{n+1})$-doubling and \rf{cc2}, we
get
\begin{eqnarray*}
\sum_j \int_{R_{s_j}} |\alpha_{s_j}^N|\,d\mu & \leq & \sum_j
\int |f\,w_{s_j}|\,d\mu\\
& \leq & C \int_{3R_k} |f|\,d\mu
\, \leq \, C \lambda \mu(6R_k)
\, \leq \, C_{14}\lambda\,\mu(R_k).
\end{eqnarray*}
Therefore,
$$\mu\left\{{\textstyle \sum_j} |\alpha_{s_j}^N| > 2C_{14}\lambda\right\}\leq
\frac{\mu(R_k)}{2}.$$ So we set
$$A_k^N = R_k\cap\left\{{\textstyle \sum_j} |\alpha_{s_j}^N| \leq
2C_{14}\lambda\right\},$$
and then $\mu(A_k^N)\geq \mu(R_k)/2.$

The constant $a_k^N$ is chosen so that for $\alpha_k^N = a_k^N\,
\chi_{A_k^N}$ we have $\int\alpha_k^N\,d\mu = \int f\,w_k\,
d\mu$.  Then we obtain
\begin{eqnarray*}  |a_k^N|  &\leq&  \frac{1}{\mu(A_k^N)}\int
|f\,w_k|\,d\mu
\leq  \frac{2}{\mu(R_k)}\int |f\,w_k|\,d\mu  \\  & \leq &
\frac{2}{\mu(R_k)}\int_{\frac{1}{2}R_k} |f|\,d\mu  \leq  C_{15}\lambda
\end{eqnarray*}
(this calculation also applies to $k=1$).
Thus,
$$|\alpha_k^N|+\sum_{j} |\alpha_{s_j}^N| \leq (2C_{14}+C_{15})\,\lambda.$$
If we choose $B=2C_{14}+C_{15}$, \rf{cc5} follows for the cubes 
$R_1,\ldots,R_n$.

Now it is easy to check that \rf{cc4.5} also holds. Indeed we have
$$\|\alpha_i^N\|_{L^\infty(\mu)} \,\mu(R_i)
\leq  C\,|a_i^N|\,\mu(A_i^N)
=  C\,\left|\int_{Q_i}f\,w_i\,d\mu\right| \leq
C\,\|\alpha_i^N\|_{L^1(\mu)}.$$

Finally, taking weak limits in the weak-$\ast$ topology of $L^\infty(\mu)$, one
easily obtains the required functions $\alpha_i$. The details are left to
reader. A similar argument can be found in the proof of Lemma 7.3 of
\cite{Tolsa3}.
\end{proof}

Using the decomposition above we can prove Lemma \ref{final} partially. This
will be the first step of its proof.

\begin{lemma}  \label{final0}
The subspace $H^1_\Phi(\mu)\cap L^\infty(\mu)$ is dense in $H^1_\Phi(\mu)$.
\end{lemma}

\begin{proof} Given $f\in H^1_\Phi(\mu)$,
for each integer $k\geq0$, we consider the generalized
Calder\'on-Zygmund decomposition of $f$ given in the preceding lemma, with
$\lambda=2^k$.
We will adopt the convention that all the elements of that decomposition will
carry the subscript $k$. Thus we write $f= g_k + b_k$, as in (c) of Lemma
\ref{final}. We know that $g_k$ is bounded and satisfies $\int g_k\,d\mu=0$
(because $\int b_k\,d\mu=0$).
We will show that $g_k\to f$ in $L^1(\mu)$ and $\|M_\Phi (g_k -
f)\|_{L^1(\mu)}\to 0$ as $k\to \infty$ too.

It is not difficult to check that $b_k$ tends to $0$ in $L^1(\mu)$.
Indeed, if we set $\Omega_k= \bigl\{M_{(2)} f(x)> 2^{k}\bigr\}$, then
$\mu(\Omega_k)\to 0 $ as $k\to\infty$, because $f\in L^1(\mu)$. Thus
$$\int |b_k|\,d\mu \leq 2 \sum_i \int |f\,w_{i,k}|\,d\mu \leq C\,
\int_{\Omega_k} |f|\,d\mu
\xrightarrow{k\to\infty}
0,$$
and so $g_k\to f$ in $L^1(\mu)$.

Let us see that $\|M_\Phi b_k\|_{L^1(\mu)}\to 0$ as $k\to \infty$. We denote
$b_{i,k} = f\,w_{i,k} - \alpha_{i,k}$. Then we have
$$\|M_\Phi b_k\|_{L^1(\mu)} \leq \sum_i  \|M_\Phi b_{i,k} \|_{L^1(\mu)}.$$
The estimates for each term $\|M_\Phi b_{i,k} \|_{L^1(\mu)}$
are (in part) similar to the ones in Lemma
\ref{only} for estimating $M_\Phi$ over atomic blocks. We write
\begin{eqnarray} \label{tyy0}
\|M_\Phi b_{i,k} \|_{L^1(\mu)} & \leq &
\int_{\R^d\setminus 2R_{i,k}} M_\Phi b_{i,k}\,d\mu \nonumber \\ && \mbox{} +
\int_{2R_{i,k}} M_\Phi (f\,w_{i,k})\,d\mu +
\int_{2R_{i,k}} M_\Phi \alpha_{i,k} \,d\mu
\end{eqnarray}
Taking into account that $\int b_{i,k}\,d\mu = 0$, it is easily seen that
$$\int_{\R^d\setminus 2R_{i,k}} M_\Phi b_{i,k}\,d\mu \leq
C\,\|b_{i,k}\|_{L^1(\mu)} \leq C\,\|f\,w_{i,k}\|_{L^1(\mu)}$$
(the calculations are similar to the ones in
\rf{tyy1} and \rf{tyy2}).

Let us consider the last term on the right hand side of \rf{tyy0} now. By
\rf{cc4} and \rf{cc4.5} we get
$$\int_{2R_{i,k}} M_\Phi \alpha_{i,k} \,d\mu \leq
 \|\alpha_{i,k}\|_{L^\infty(\mu)} \mu(2R_{i,k})\,d\mu
\leq C\,\|f\,w_{i,k}\|_{L^1(\mu)}.$$

We split the second integral on the right hand side of \rf{tyy0} as follows:
$$\int_{2R_{i,k}} M_\Phi (f\,w_{i,k})\,d\mu = \int_{2R_{i,k}\setminus
2Q_{i,k}} + \int_{2Q_{i,k}}.$$
As in \rf{tyy3}, we have
\begin{eqnarray*}
\int_{2R_{i,k}\setminus 2Q_{i,k}} M_\Phi (f\,w_{i,k})\, d\mu & \leq &
C\,\|f\,w_{i,k}\|_{L^1(\mu)}\, \int_{2R_{i,k}\setminus 2Q_{i,k}}
\frac{1}{|x-z_{Q_{i,k}}|^n}\, d\mu(x)  \\
& \leq & C\,\|f\,w_{i,k}\|_{L^1(\mu)}
\, (1+ \delta(Q_{i,k},R_{i,k})) \\
& \leq & C\,\|f\,w_{i,k}\|_{L^1(\mu)}.
\end{eqnarray*}

Finally we have to deal with $\int_{2Q_{i,k}}
M_\Phi (f\,w_{i,k})\, d\mu$. Consider $x\in 2Q_{i,k}$ and $\vphi\sim x$. Then
\begin{equation} \label{tyy4}
\left|\int \vphi\, (f\,w_{i,k}) \,d\mu \right| =
\left|\int (\vphi\, w_{i,k})\,f \,d\mu \right| \leq C\,M_\Phi f(x),
\end{equation}
because $C\,\vphi\,w_{i,k}\sim x$ for some constant $C>0$. Indeed, for
$y\in\R^d$ we have
$$0\leq w_{i,k}\vphi(y) \leq \vphi(y) \leq \frac{1}{|y-x|^n}$$
and
\begin{eqnarray*}
|(\vphi\,w_{i,k})'(y)| & \leq & |\vphi'(y)\,w_{i,k}(y)| +
|\vphi(y)\,w_{i,k}'(y)| \\
& \leq & \frac{1}{|y-x|^{n+1}} + \frac{C}{|y-x|^{n}}\,|w_{i,k}'(y)|.
\end{eqnarray*}
Recall also that $|w_{i,k}'(y)| \leq C\ell(Q_{i,k})^{-1}$ and $\supp(w_{i,k})
\subset 2Q_{i,k}$. Then we get $|w_{i,k}'(y)| \leq C\,|y-x|^{-1}$ for all
$y\in\R^d$. Thus $|(\vphi\,w_{i,k})'(y)| \leq C\,|y-x|^{-n-1}$.
So \rf{tyy4} holds and then
$$\int_{2Q_{i,k}}
M_\Phi (f\,w_{i,k})\, d\mu \leq C\,\int_{2Q_{i,k}} M_\Phi f\, d\mu.$$

When we gather the previous estimates, we obtain
$$\|M_\Phi b_{i,k} \|_{L^1(\mu)} \leq C\,\|f\,w_{i,k}\|_{L^1(\mu)} +
C\,\int_{2Q_{i,k}} M_\Phi f\, d\mu.$$
Taking into account the finite overlap of the cubes $2Q_{i,k}$ (recall that
they are Whitney cubes covering $\Omega_k$), we get
$$\|M_\Phi b_k\|_{L^1(\mu)} \leq C\, \int_{\Omega_k} (|f| + M_\Phi f)\, d\mu
\xrightarrow{k\to\infty} 0,$$
and we are done.
\end{proof}


\vspace{3mm}
\noindent{\bf Proof of Lemma \ref{final}}.
Take $f\in H^1_\Phi(\mu)\cap L^\infty(\mu)$.
Consider the infinite increasing sequence of the cubes 
$Q_k = 4^{N_k}\,[-1,1]^d$ that are $(4,4^{n+1})$-doubling. Let $w$ be a 
$\CC^\infty$ function such that $\chi_{[-1,1]^d}(x) \leq 
w(x) \leq \chi_{[-2,2]^d}(x)$ for all $x$. We denote $w_k(x) = w(4^{-N_k}x)$
(so $\chi_{Q_k}(x) \leq w_k(x)\leq \chi_{2Q_k}(x)$) and we set 
$$f_k = w_k\,f - \frac{\chi_{Q_k}}{\mu(Q_k)}\,\int w_k\,f\,d\mu.$$
It is clear that $f_k$ is bounded, has compact support and converges to $f$ in
$L^1(\mu)$ as $k\to\infty$. We will prove that
\begin{eqnarray} \label{acopri}
\|M_\Phi (f-f_k)\|_{L^1(\mu)} & \leq & C\,\left|\int w_k\,f\,d\mu\right|
+ C\, \int_{\R^d\setminus 4Q_k} M_\Phi f\,d\mu  \\
&&\mbox{} + \int_{4Q_k} M_\Phi
((1-w_k)\,f)\, d\mu. \nonumber
\end{eqnarray} 
Finally we will show that the terms on the right hand side of \rf{acopri} tend
to $0$ as $k\to\infty$ and we will be done.

Let us consider first the integral of $M_\Phi (f-f_k)$ over $\R^d\setminus
4Q_k$. We set
$$\int_{\R^d\setminus 4Q_k} M_\Phi (f-f_k)\,d\mu \leq 
\int_{\R^d\setminus 4Q_k} M_\Phi f\,d\mu + 
\int_{\R^d\setminus 4Q_k} M_\Phi f_k\,d\mu.$$
We only have to estimate the last integral on the right hand side.
Take $x\in\R^d\setminus 4Q_k$, $\vphi\sim x$ and let
$y_0\in 2Q_k$ be the point where $\vphi$ attains its minimum over $2Q_k$ (recall
that we assume $\vphi\geq0$ and $\vphi\in\CC^1$). We denote $c_k = 
\int w_k\,f\,d\mu / \mu(Q_k)$ and then we set
\begin{eqnarray*}
\int f_k\,\vphi\,d\mu & = & \int f(y) \,(\vphi(y)-\vphi(y_0))\,d\mu(y) \\
& = & \int w_k(y)\,f(y)\,(\vphi(y)-\vphi(y_0))\,d\mu(y) \nonumber \\
&&\mbox{}- c_k \int_{Q_k}(\vphi(y)-\vphi(y_0))\,d\mu(y) \,= \,I_1- I_2.
\end{eqnarray*}
Let us consider the function $\psi(y) = w_k(y)\,(\vphi(y) - \vphi(y_0))$. This
function satisfies
$$0\leq \psi(y) \leq \vphi(y)$$
and
\begin{eqnarray*}
|\psi'(y)| & \leq & |w_k(y)\,\vphi'(y)| + |w_k'(y)|\,|\vphi(y) - \vphi(y_0)| \\
& \leq & \frac{1}{|y-x|^{n+1}} + 
C\,\ell(Q_k)^{-1}\, \frac{\ell(Q_k)}{|y-x|^{n+1}} \,=\,C\,
\frac{1}{|y-x|^{n+1}}.
\end{eqnarray*}
Therefore $C\,\psi\sim x$ for some constant $C>0$ and so
$|I_1| \leq C\,M_\Phi f(x).$
For $I_2$ we use a cruder estimate:
$$|I_2| \leq C\,|c_k|\,\mu(Q_k)\,\frac{\ell(Q_k)}{|y_0-x|^{n+1}}.$$
Thus we obtain
$$M_\Phi f_k(x) \leq C\,M_\Phi f(x) + 
C\,|c_k|\,\mu(Q_k)\,\frac{\ell(Q_k)}{|y_0-x|^{n+1}}.$$
Since
$$\int_{\R^d\setminus 4Q_k} \frac{1}{|y_0-x|^{n+1}}\,d\mu(x)
\leq C\,\ell(Q_k)^{-1},$$
we get
\begin{eqnarray}  \label{tgbb0}
\int_{\R^d\setminus 4Q_k} M_\Phi f_k\,d\mu & \leq &
C\, \int_{\R^d\setminus 4Q_k} M_\Phi f\,d\mu + C\,|c_k|\,\mu(Q_k)\nonumber \\
& = & C\, \int_{\R^d\setminus 4Q_k} M_\Phi f\,d\mu + C\,
\left|\int w_k\,f\,d\mu\right|.
\end{eqnarray}

Now we have to deal with $\int_{4Q_k} M_\Phi (f-f_k)\,d\mu$.
For $x\in 4Q_k$ we write
\begin{equation} \label{tgbb1}
M_\Phi (f-f_k)(x) \leq M_\Phi ((1-w_k)\,f)(x) + \,
M_\Phi \left(\frac{|c_k|}{\mu(Q_k)}\,\chi_{Q_k}\right)(x).
\end{equation}
Since $M_\Phi \chi_{Q_k}(x)\leq 1$ and $Q_k$ is $(4,4^{n+1})$-doubling, we get
\begin{equation}  \label{tgbb2}
\int_{4Q_k} M_\Phi \left(\frac{|c_k|}{\mu(Q_k)}\,\chi_{Q_k}\right)(x)\,d\mu(x)
\leq C\,|c_k| = C\,\left|\int w_k\,f\,d\mu\right|.
\end{equation}
From \rf{tgbb0}, \rf{tgbb1} and \rf{tgbb2} we derive \rf{acopri}.

Now we have to see that the terms on the right hand side of \rf{acopri} tend to $0$
as $k\to\infty$. Since $f,\,M_\Phi f\in L^1(\mu)$, by the dominated convergence
theorem
$$\lim_{k\to\infty}\left|\int w_k\,f\,d\mu\right| +
\int_{\R^d\setminus 4Q_k} M_\Phi f\,d\mu = 0.$$
Let us turn our attention to the third term on the right hand side of
\rf{acopri}. Take $x\in 4Q_k$ and $\vphi\sim x$.
It is easily seen that $C\,w_k\,\vphi\sim x$ for some
constant $C>0$. So we get $M_\Phi (w_k\,f)(x) \leq C\,M_\Phi f(x)$ and then for
any $x\in \R^d$,
$$\chi_{4Q_k}(x)\,M_\Phi ((1-w_k)\,f)(x) \leq \chi_{4Q_k}(x)\,
(M_\Phi f(x) + M_\Phi (w_k\,f)(x)) \leq C\,M_\Phi f(x).$$
Therefore, if we show that $\chi_{4Q_k}(x)\,M_\Phi ((1-w_k(x))\,f)(x)$ tends to
$0$ pointwise as $k\to\infty$, we will be done by a new application of the
dominated convergence theorem.

For a fixed $x\in\R^d$, let $k_0$ be such that $x\in \frac{1}{2}Q_k$ for $k\geq
k_0$. Notice that if $\vphi\sim x$ and $y\not\in Q_k$, 
then $|\vphi(y)|\leq C/\ell(Q_k)^n$. Thus
\begin{eqnarray*}
\left|\int \vphi(x)(1-w_k(x))\,f(x)\,d\mu(x)\right| & \leq &
\|f\|_{L^1(\mu)}\,\|(1-w_k)\,\vphi\|_{L^\infty(\mu)} \\
& \leq & C\,\frac{\|f\|_{L^1(\mu)}}{\ell(Q_k)^n}.
\end{eqnarray*}
Then we get 
$$\chi_{4Q_k}(x)\,M_\Phi ((1-w_k(x))\,f)(x) \leq 
C\,\frac{\|f\|_{L^1(\mu)}}{\ell(Q_k)^n} \xrightarrow{k\to\infty} 0.$$
\qed




\vspace{10mm}


\begin{thebibliography}{MMNO}

\bibitem[Ca]{Carleson} L. Carleson. {\em Two remarks on $H^1$ and $BMO$,}
Advances in Math. {\bf 22} (1976), 269-277.

\bibitem[Co]{Coifman} R.R. Coifman. {\em A real variable characterization
of $H^p$,} Studia Math. {\bf51} (1974) 269-274.

\bibitem[CW]{CW} R.R. Coifman, G. Weiss. {\em Extensions of Hardy spaces and
their use in analysis,} Bull. Amer. Math. Soc. {\bf83} (1977), 569-645.

\bibitem[GR]{GR} J. Garc\'{\i}a-Cuerva, J.L. Rubio de Francia.
{\em Weighted norm
inequalities and related topics,} North-Holland Math. Studies 116, 1985.

\bibitem[Jo]{Journe} J.-L. Journ\'e. {\em Calder\'on-Zygmund operators,
pseudo-differential operators and the Cauchy integral of Calder\'on.}
Lecture Notes in Math. 994, Springer-Verlag, 1983.

\bibitem[La]{Latter} R.H. Latter. {\em A characterization of $H^p(\R^n)$ in
terms of atoms,} Studia Math. {\bf62} 1978, 92-101.

\bibitem[MS1]{MS1} R.A. Mac\'{\i}as, C. Segovia.
{\em Lipschitz functions on spacesof homogeneous type,} Advances in Math.
{\bf33} (1979), 257-270.

\bibitem[MS2]{MS2} R.A. Mac\'{\i}as, C. Segovia. {\em A decomposition into atoms
of distributions on spaces of homogeneous type,} Advances in Math. {\bf33}
(1979), 271-309.

\bibitem[MMNO]{MMNO} J. Mateu, P. Mattila, A. Nicolau, J. Orobitg. {\em
$\bmo$ for non doubling measures.} To appear in Duke Math. J.

\bibitem[NTV1]{NTV1}  F. Nazarov, S. Treil, A. Volberg. {\em Cauchy integral
and Calder\'on-Zygmund operators on nonhomogeneous spaces,}
Int. Math. Res. Not. {\bf 15} (1997), 703-726.

\bibitem[NTV2]{NTV2} {F. Nazarov, S. Treil, A. Volberg.}
{\em Weak type estimates and Cotlar inequalities for
Calder\'on-Zygmund operators in nonhomogeneous spaces,}
Int. Math. Res. Not. {\bf 9} (1998), 463-487.

\bibitem[NTV3]{NTV*}
F. Nazarov, S. Treil, A. Volberg. {\em $Tb$-theorem on non-homogeneous
spaces.} Preprint (1999).

\bibitem[NTV4]{NTV4} F. Nazarov, S. Treil, A. Volberg. {\em Accretive
$Tb$-systems on non-homogeneous spaces.} Preprint (1999).

\bibitem[OP]{OP} J. Orobitg, C. P\'erez. {\em $A_p$ weights for non doubling
measures in $\R^n$ and applications.} Preprint (1999).

\bibitem[St]{Stein} E. M. Stein. {\em Harmonic analysis. Real-Variable
methods, orthogonality, and oscillatory integrals.} Princeton Univ. Press.
Princeton, N.J., 1993.

\bibitem[To1]{Tolsa1} X. Tolsa. {\em $L^2$-boundedness of the Cauchy integral
operator for continuous measures,} Duke Math. J. {\bf 98}:2 (1999), 269-304.


\bibitem[To2]{Tolsa2.5} X. Tolsa. {\em A $T(1)$ theorem for non doubling
measures with atoms.} Proc. London Math. Soc. To appear.

\bibitem[To3]{Tolsa3} X. Tolsa. {\em $\bmo$, $H^1$ and Calder\'on-Zygmund
operators for non doubling measures}. Preprint (1999). Provisionally accepted
for publicaton in Math. Ann.

\bibitem[To4]{Tolsa5} X. Tolsa. {\em A proof of the weak $(1,1)$ inequality
for singular integrals with non doubling measures based on a Calder\'on-Zygmund
decomposition.} Preprint (1999).

\bibitem[Ve]{Verdera} J. Verdera, {\em On the $T(1)$ theorem for the
Cauchy integral,} Arkiv f. Mat. To appear.

\bibitem[Uc]{Uchiyama} A. Uchiyama. {\em A maximal function characterization
of $H^p$ on the space of homogeneous type}, Trans. Amer. Math. Soc. {\bf262}:2,
(1980), 579-592.

\end{thebibliography}
\end{document}